\documentclass[english,12pt,twoside]{article}
\usepackage[a4paper, margin=2.1cm]{geometry}
\usepackage[english]{babel}

\usepackage{mathrsfs}
\usepackage[ruled,vlined]{algorithm2e}
\usepackage{color}
\usepackage{mathtools}

\usepackage{enumitem}

\usepackage{amsmath}
\usepackage{amsfonts}
\usepackage{amssymb}

\usepackage{hyperref}
\usepackage{relsize}

\usepackage{soul}
\usepackage{apacite}
\usepackage{comment}

\usepackage{amsthm}
\newtheorem{theorem}{Theorem}[]
\newtheorem{lemma}{Lemma}[]
\newtheorem{corollary}{Corollary}[]
\theoremstyle{definition}
\newtheorem{definition}{Definition}[]

\usepackage{thmtools}
\declaretheorem[style=definition,qed=$\triangle$]{example}

\newcommand{\argmax}{\operatornamewithlimits{argmax}}
\newcommand{\argmin}{\operatornamewithlimits{argmin}}

\usepackage{tikz}
\usetikzlibrary{arrows.meta}
\usepackage{graphics}
\usepackage{subcaption}
\usepackage[round]{natbib}
\begin{document}


\title{Finding optimal Stackelberg production strategies: \\ How to produce in times of war?
}
\author{Loe Schlicher\\
\footnotesize{School of Industrial Engineering, Eindhoven University of Technology} \\
\footnotesize{P.O. Box 513, 5600 MB, Eindhoven, The Netherlands} \\
\footnotesize{l.p.j.schlicher@tue.nl}
\and
Marieke Musegaas\\
  \footnotesize{
  Department of Data Science and Knowledge Engineering, Maastricht University} \\
  \footnotesize{P.O. Box 616, 6200 MD, Maastricht, The Netherlands } \\
\footnotesize{m.musegaas@maastrichtuniversity.nl}
\and
Herman Blok\\
\footnotesize{Department of Mathematics, Vrije Universiteit Amsterdam} \\
\footnotesize{De Boelelaan 1081a, 1081 HV,  Amsterdam, The Netherlands} \\
\footnotesize{hermanblok@gmail.com}
}

\date{\today}

\maketitle


\begin{abstract}
\noindent Inspired by a military context, we study a Stackelberg production game  where a country's government, the leader, wants to maximize the production of military assets. The leader does so by allocating his resources among a set of production facilities. His opponent, the follower, observes this allocation and tries to destroy the associated production as much as possible by allocating his destructive resources, for example bombs, among these facilities. In this paper, we identify a follower's optimal strategy. For the leader, we show that an optimal production strategy can be found in the class of so-called seried-balanced strategies. We present a linear time algorithm that finds an optimal strategy in this class.
\end{abstract}

\section{Introduction}

Strategic bombing is a military strategy that is, among others, used to destroy military production facilities in times of war. In order to be less vulnerable to such bombings, a government can decide to spread its military production over multiple production facilities. Inspired by this military context, we investigate in this paper how a government should optimally spread the  production of a single type of military asset over its production facilities, while taking into account the impact of strategic bombing. We do so by introducing and solving a Stackelberg production game that models such a military setting in times of war.

In this Stackelberg production game, a country's government (the leader) has to allocate a fixed amount of resources among several production facilities, in order to maximize the production of a single type of military asset (e.g., tanks). For this, the leader takes into account that each production facility has its own production rate. That means,  some production facilities can produce more assets per resource than others. Once the leader has allocated his resources, the country's opponent (the follower) observes this allocation and tries to destroy the leader's production as much as possible. The follower does so by allocating a fixed amount of destructive resources (e.g., bombs) among the production facilities, which will destroy the production facilities and produced military assets. For this, the follower takes into account that some production facilities are harder to destroy than others.


In this paper, we study the leader's and follower's optimal strategies. First, we identify an intuitive optimal strategy for the follower. We show that in this optimal strategy, the follower destroys those production facilities that result in the highest production loss per destructive resource. Thereafter, we present a linear time algorithm to find an optimal production strategy for the leader. This algorithm, which is the main contribution of this paper, selects a production strategy from the  class of so-called seried-balanced strategies. We prove the correctness of this algorithm by showing that there always exists an optimal production strategy in this class.


As the name suggests, the production strategies in this class satisfy two properties. The first property, balancedness, states that all resources are allocated among a subset of production facilities such that it is equally likely for the follower to destroy any of these. The second property, seriedness, states that this subset consists of those production facilities with the highest production rates. Balancedness is a commonly seen property in, among others, military-oriented games (see, e.g., \cite{powell2009sequential}). However, seriedness is not, and may be counter intuitive. In particular, this property states that it is always beneficial from the leader’s
perspective to invest some resources in a production facility with the highest production rate. This is remarkable, since there is no direct relationship between the difficulty of destroying this production facility and its production rate. In this paper, we explain why it is optimal to select a production strategy from the class of seried-balanced strategies.

In the literature, several military-oriented Stackelberg games have been studied. In these games, the leader and follower typically optimize opposed objective functions (e.g., maximizing/minimizing military production as in this paper). Since  the follower chooses the solution which leads to the worst case objective function value for the leader, these military-oriented games are  typically expressed as bilevel optimization problems (see, e.g., \cite{dempe2002}).

Military-oriented Stackelberg games have shown to be relevant for several types of settings. For instance, \cite{israeli2002shortest} study how to route a military supply convoy through a hostile environment and \cite{morton2007models} illustrate how to optimally protect a border against smuggling. In another setting, \cite{powell2009sequential} and \cite{scaparra2008exact} study how a government must protect its vulnerable locations (e.g., buildings or train stations) against an opponent that plans to destroy these locations. As a final example,  \cite{gutin2014interdiction} show how a government can delay an opponent's nuclear weapons development project.






The literature on military-oriented Stackelberg games can be divided in two main streams. The first one focuses on finding techniques or methods to solve the underlying bilevel optimization problems efficiently (see, e.g., \citet{israeli2002shortest,scaparra2008bilevel,scaparra2008exact,liberatore2011analysis, cappanera2011optimal}; \citet{gutin2014interdiction,jiang2018multi,starita2016optimizing,washburn1995two}).
Complementary to these computational-oriented studies,
there also exists a stream of literature that focuses on the analytical tractability of Stackelberg games (see, e.g., \citet{bier2007choosing, powell2009sequential,zhuang2011secrecy,hausken2011defending,shan2013hybrid}). In these studies, the focus is on the identification and characterization of the leader’s and follower’s optimal strategies. In this paper,  we also identify  optimal  strategies for both the leader and the follower,  and  thus we contribute  to  this  second  stream of literature.

Our Stackelberg production game, and the associated results, may also find applications in other fields. For instance, in a retail setting, our game illustrates how to distribute products across countries, while taking into account the possibility of supply disruptions. In such a case, the retailer could be recognized as the leader, the production rates as the profit margins per country and the follower's decision as the worst case supply disruption.

The organization of this paper is as follows. Section~\ref{section:Thestackelbergproductiongame} formally introduces the Stackelberg production game. In Section~\ref{sectie:The follower's best replies} we focus on the identification of a follower's optimal strategy and in Section~\ref{sectie: The leader's best replies} we focus on the identification of a leader's optimal production strategy. Finally, in Section~\ref{sectie: conclusions} we summarize our results and provide three directions for future research.

\section{The Stackelberg production game} \label{section:Thestackelbergproductiongame}


We consider a Stackelberg production game where a country's government is the \emph{leader} and the country's opponent is the \emph{follower}. The leader has to allocate $R_l \in \mathbb{R}_{> 0}$ resources
among a set of $N=\{1,2,\ldots,n\}$ facilities, in order to maximize the production of a single type  of military asset (e.g., tanks). Once the leader has allocated his resources, the follower observes this allocation and tries to destroy the produced military assets of the leader as much as possible. The follower does so by allocating $R_f \in \mathbb{R}_{> 0}$ destructive resources (e.g., bombs) among the facilities, which will destroy the facilities and thus also the military assets produced in this facility.

A \emph{strategy for the leader} is defined as an allocation of his resources to the facilities. We denote this strategy by $x = (x_i)_{i \in N} \in \mathbb{R}^n_{\geq 0}$. Since the leader has to allocate
$R_l$ resources, such a strategy needs to satisfy $\sum_{i \in N} x_i \leq R_l$. We denote the set of all feasible strategies by
$$\mathscr{X} = \left\{x \in \mathbb{R}^n_{\geq 0}~\Bigg\vert~ \sum_{i=1}^n\ x_i \leq R_l\right\}.$$
We assume that each facility has a unique production rate. This \emph{production rate}, denoted by $p_i \in \mathbb{R}_{>0}$, is the number of military assets produced per invested resource. Hence, if the leader allocates $x_i$ resources to facility $i$, a \emph{facility production} $p_i x_i$ is obtained.

A \emph{strategy for the follower} is defined as an allocation of his destructive resources to the facilities. We denote this strategy by $y = (y_i)_{i \in N} \in \mathbb{R}^n_{\geq 0}$.  Since the follower has only $R_f$ destructive resources, such a strategy needs to satisfy  $\sum_{i \in N} y_i \leq R_f$. Moreover, we assume that each facility has a unique destruction quantity. This \emph{destruction quantity}, denoted by  $a_i \in \mathbb{R}_{>0}$, is the amount of destructive resources needed in order to completely destroy a facility and its associated production. If the follower allocates $y_i$ destructive resources to facility $i$,  a \emph{facility production reduction} of $p_ix_i \frac{y_i}{a_i} $ is obtained. We assume that the follower never allocates more destructive resources than needed. Hence, in addition to the resource constraint, we also require $y_i \leq a_i$. We denote the set of all feasible strategies by
\begin{equation}
\label{eq:feasiblestategiesfollower}
    \mathscr{Y} = \left\{y \in \mathbb{R}^n_{\geq 0}~\Bigg\vert~\sum_{i=1}^n\ y_i \leq R_f, y_i \leq a_i \text{ for all } i \in N \right\}.
\end{equation}
 In order to exclude the trivial case, where the follower has enough destructive resources to destroy all facilities of the leader, we assume $R_f < \sum_{i \in N}{a_i}$.

For a given strategy $x \in \mathscr{X}$ of the leader and a given strategy  $y \in \mathscr{Y}$ of the follower, the associated leader's \emph{total production after destruction} is given by the sum of the \textcolor{black}{facility productions} minus the sum of the \textcolor{black}{facility production reductions}, i.e.,
\begin{equation}
\label{vgl:Loe}
\mathscr{P}(x,y)
=\sum_{i \in N}{ p_ix_i} - \sum_{i \in N}{p_ix_i \tfrac{y_i}{a_i} }.
\end{equation}
The following example illustrates how to calculate \textcolor{black}{facility productions}, \textcolor{black}{facility production reductions} and total production after destruction. The Stackelberg production game in this example is used throughout the entire paper.

\begin{example} \label{example1} Consider the Stackelberg production game  with $N = \{1,2,3,4,5\}$, $R_l=5$ and $R_f=1 \frac{3}{4}$. The production rates  $\{p_i\}_{i \in N}$  and destruction quantities  $\{a_i\}_{i \in N}$ are presented in Table \ref{table1a}.

\begin{table}[ht]  \centering
\begin{tabular}{lccccc} \hline
Facility $i$ & 1 & 2 & 3 & 4 & 5 \\ \hline
$p_i$ & 12 & 8 & 5 & 2 & 1 \\
$a_i$ & $\frac{9}{10}$  & 1 & $\frac{1}{4}$ &  1 & $\frac{3}{4}$ \\
\hline
 \end{tabular}
 \caption{The production rates and destruction quantities of Example~\ref{example1}.}
 \label{table1a}
 \end{table}

Suppose the leader applies strategy $x= \left(0,\frac{7}{10},\frac{3}{10},0,4\right)$ and the follower applies strategy $y = \left(0,\frac{7}{8},\frac{1}{8}, 0,\frac{3}{4}\right)$. Note that indeed both strategies are feasible, i.e., $x \in \mathscr{X}$ and $y \in \mathscr{Y}$. The resulting facility productions $\{p_ix_i\}_{i \in N}$ and facility production reductions $\{p_ix_i\frac{y_i}{a_i}\}_{i \in N}$ are as presented in Table \ref{table1b}.

\begin{table}[ht]  \centering
\begin{tabular}{lccccc} \hline
Facility $i$ & 1 & 2 & 3 & 4 & 5 \\ \hline
$p_ix_i$& 0 &  $5 \frac{6}{10}$ & $1\frac{1}{2}$ & 0 & $4 $  \\
$p_ix_i\frac{y_i}{a_i} $ & 0 & $4 \frac{9}{10}$ & $\frac{3}{4}$ & 0 & $4$ \\\hline
 \end{tabular}
 \caption{Facility productions and facility production reductions of Example~\ref{example1}.}
 \label{table1b}
 \end{table}

\noindent Using Table~\ref{table1b} we conclude that the total production after destruction equals
\begin{equation*}
\mathscr{P}(x,y) = \left(5 \tfrac{6}{10} + 1\tfrac{1}{2} + 4\right) - \left( 4 \tfrac{9}{10} + \tfrac{3}{4} + 4\right) = 1 \tfrac{9}{20}.~\qedhere
\end{equation*}
\end{example}

The follower has the aim to minimize the leader's total production after destruction. We denote this resulting \emph{worst case total production after destruction}, for a given strategy $x \in \mathscr{X}$ of the leader, by
$$\mathscr{P}(x)=\min_{y \in \mathscr{Y}}{\mathscr{P}(x,y)}.$$
The leader has the possibility to anticipate on this outcome. In particular, the leader will anticipate by maximizing the worst case total production after destruction, which is formally defined as
$$\max_{x \in \mathscr{X}}{\mathscr{P}(x)} = \max_{x \in \mathscr{X}}{\min_{y \in \mathscr{Y}}{\mathscr{P}(x,y)}}.$$

In this paper, we focus on the identification of an optimal strategy for the leader. Given the sequential order in which the players move in this Stackelberg production game, we start in Section~\ref{sectie:The follower's best replies}  with identifying an optimal strategy for the follower. Thereafter, in Section~\ref{sectie: The leader's best replies}, we  present a linear time algorithm to find an optimal strategy for the leader. This algorithm, which is the main contribution of this paper, selects a strategy from  the  class  of  seried-balanced strategies.

\section{Optimal strategy for the follower}
\label{sectie:The follower's best replies}


\noindent In this section we identify an optimal strategy for the follower. We show that this optimal strategy has a very easy and intuitive structure.

Recall that the set of all feasible strategies for the follower, denoted by $\mathscr{Y}$, is such that in total at most $R_f$ destructive resources are allocated and moreover at most $a_i$ destructive  resources are allocated to facility $i \in N$ (cf.~\eqref{eq:feasiblestategiesfollower}).  Given a strategy $x \in \mathscr{X}$ for the leader, the follower wants to allocate the destructive resources in such a feasible way that the leader's total production after destruction is as low as possible. It follows from~\eqref{vgl:Loe} that this is equivalent to maximizing $\sum_{i \in N}{p_ix_i\frac{y_i}{a_i}}$.  In Theorem~\ref{thm:optimalstrategyfollower} we state  that it is optimal for the follower to destroy the facility with the highest ratio $\frac{p_ix_i}{a_i}$ first, thereafter destroy the facility with the second highest ratio $\frac{p_ix_i}{a_i}$, and so on, until no destructive resources are left. Note that, due to ties in these ratios, there might be different non-increasing orders and thus different optimal strategies. We now formally describe this set of optimal strategies.

Given a leader's strategy $x \in \mathscr{X}$, we define the \emph{destruction ratio} as the facility's production reduction per destructive resource. Hence, the destruction ratio for facility $i \in N$ is given by the ratio $\frac{p_ix_i}{a_i}$. Next, let $\Omega^x$ denote the set of all permutations of $N$ that orders the facilities based on this destruction ratio in non-increasing order, i.e.,
\begin{equation*} 
\Omega^x = \left\{ \omega \in \Omega \hspace{1mm} \bigg \vert \hspace{1mm} \frac{p_{\omega(i)}x_{\omega(i)}}{a_{\omega(i)}} \geq \frac{p_{\omega(j)} x_{\omega(j)}}{a_{\omega(j)}} \mbox{ for all } i,j \in N \mbox{ with } i<j\right\}, \end{equation*}
where $\Omega$ denotes the set of all permutations of $N$. Then, every permutation $\omega \in \Omega^x$ corresponds to an optimal strategy for the leader. Namely, allocate $a_{\omega(1)}$ destructive resources to facility $\omega(1)$ first, then allocate $a_{\omega(2)}$ destructive resources to facility $\omega(2)$, and so on, until we reach the last facility for which there are still destructive resources available. We denote this facility by $q^{\omega}$, i.e.,
\begin{equation*}
q^{\omega} = \min\left\{j \in N \hspace{1mm} \bigg \vert \hspace{1mm} \sum_{i=1}^j a_{\omega(i)} \geq R_f  \right\}.
\end{equation*}
Note that facility $q^{\omega}$ always  exists due to the assumption $R_f < \sum_{i \in N}{a_i}$. Since this facility gets allocated the remaining amount of destructive resources, this facility is the only one that might be partly destroyed.
We denote the set of facilities that are (either fully or partly) destroyed  by $A^{\omega}$, i.e.,

\begin{equation*}
A^{\omega} = \left\{i \in N \hspace{1mm} \bigg\vert \hspace{1mm} \omega(i) \leq \omega( q^{\omega} ) \right\}.
\end{equation*}
We are now ready to formally define the optimal strategy corresponding to permutation $\omega$, denoted by $y ^{\omega}$. Namely,
\begin{equation}
\label{eq:defyomega}
y^{\omega}_i = \begin{cases}
 a_i & \text{ if } i \in A^{\omega} \backslash \{q^{\omega}\}, \\
{\displaystyle R_f-  \sum_{j \in A^{\omega} \backslash \{q^{\omega}\}}{a_j}} & \text{ if } i=q^{\omega}, \\
0 & \text{ otherwise}. \end{cases}
\end{equation}
The following theorem states that  $y^{\omega}$ is a  follower's optimal strategy. We want to emphasize that all proofs are relegated to the appendix.


\begin{theorem}
\label{thm:optimalstrategyfollower}
Let  $x \in \mathscr{X}$. Then, every strategy  $y^{\omega}$, with $\omega \in \Omega^x$, is optimal for the follower. \end{theorem}

It immediately follows from Theorem~\ref{thm:optimalstrategyfollower} that the worst case total production after destruction is, for every $ x \in \mathscr{X}$, given by
$$\mathscr{P}(x) = \mathscr{P}(x,y^{\omega}),$$
where 
$\omega \in \Omega^x$. The following example illustrates the result of Theorem~\ref{thm:optimalstrategyfollower}. 

\newpage
\begin{example} \label{example2} Assume the leader applies strategy $x= \left(0,\frac{7}{10},\frac{3}{10},0,4\right)$. The resulting destruction ratios $\{\frac{p_ix_i}{a_i}\}_{i \in N}$ are presented in Table \ref{table2}.

\begin{table}[ht] \centering
\begin{tabular}{lccccc} \hline
Facility $i$  & 1 & 2 & 3 & 4 & 5 \\ \hline
$\frac{p_ix_i}{a_i} $ & 0 & $5 \frac{6}{10}$ & 6 &  0 & $5 \frac{1}{3}$ \\\hline
\end{tabular}
\caption{The destruction ratios of Example~\ref{example2}}
\label{table2}
\end{table}

\noindent Sorting the facilities in non-increasing order based on their destruction ratios results in  $\Omega^x = \{(3,2,5,1,4), (3,2,5,4,1)\}$. Based on the destruction quantities $\{a_i\}_{i \in N}$ we can conclude that it is optimal for the follower to fully destroy facilities 2 and 3, and to destroy facility 5 with his leftover destructive resources. Hence, for every  $\omega \in \Omega^x$,  we have $q^{\omega} = 5$, $A^{\omega}=\{2,3,5\}$ and
$y^{\omega} = (0,1,\frac{1}{4},0,\frac{1}{2})$. By Theorem~\ref{thm:optimalstrategyfollower}, it follows that the worst case total production after destruction for strategy $x$ is given by
$$\mathscr{P}(x) = \mathscr{P}(x,y^{\omega}) = \sum_{i \in N}{ p_ix_i} - \sum_{i \in N}{p_ix_i \tfrac{y^{\omega}_i}{a_i} } =  (5 \tfrac{6}{10} + 1 \tfrac{1}{2} + 4) - (5 \tfrac{6}{10} +1 \tfrac{1}{2} + 2 \tfrac{2}{3}) = 1\tfrac{1}{3}.$$
Note that strategy $y$ from Example~\ref{example1} is not an optimal strategy for the follower because $\mathscr{P}(x,y) = 1 \frac{9}{20} >  1\tfrac{1}{3}  =  \mathscr{P}(x)$.
  \end{example}

\section{Optimal strategy for the leader}
\label{sectie: The leader's best replies}
In this section we develop a linear time algorithm that finds  an optimal strategy for the leader. Specifically, this algorithm selects an optimal strategy from a set of potential optimal strategies, called seried-balanced strategies. Among this set of seried-balanced strategies, the algorithm selects the one for which the worst case total production after destruction is the highest. The algorithm does so by selecting the seried-balanced strategy with the highest composed net production rate.

In order to prove the correctness of the algorithm (as formalized in  Algorithm 1 in Section~\ref{sectie:algorithm}), we first need to prove that it suffices to consider only the class of seried-balanced strategies. This proof consists of the following three steps.

\begin{itemize}
    \item We first show that each strategy is dominated by an associated  balanced or semi-balanced strategy. As a consequence, it suffices  to consider only the class of balanced and semi-balanced strategies (Section~\ref{sectie:(semi)-bal}).
    \item  Then, we show that each semi-balanced strategy is dominated by an associated balanced strategy. As a consequence, by combining this with the result from the previous step,  it suffices to consider only the class of balanced strategies (Section~\ref{sectie:bal}).
    \item  Finally, we show that each balanced strategy is dominated by an associated seried-balanced strategy. As a consequence, by combining this with the result from the previous step, it suffices to consider only the class of seried-balanced strategies (Section~\ref{sectie:seried}).
\end{itemize}

Using these results, we complete the proof of correctness of Algorithm 1 in Section~\ref{sectie:algorithm} by showing a relationship between the worst case total production after destruction and the composed net production rate. As a result, an optimal strategy can be found by selecting the seried-balanced strategy with the highest composed net production rate.

\subsection{Reducing to the class of balanced and semi-balanced strategies}
\label{sectie:(semi)-bal}

In this section we show that, for finding an optimal  strategy, it suffices to consider the class of so-called balanced and semi-balanced strategies. A strategy is  balanced if all resources are allocated among a subset of facilities such that each of them is equally likely being destroyed by the follower. That means, the destruction ratios of the facilities in this subset are the same. Definition 1 provides the formal definition of a balanced strategy.

 \begin{definition}
 A strategy $x$ is  \emph{balanced} if there exists a set of facilities $S^x \subseteq N$  such that
 \begin{enumerate}[leftmargin=1.3cm]
\item[(i)]  $\sum_{i \in N}{x_i} = R_l$,
     \item[(ii)] $x_i>0$  for all $i \in S^x$,
     \item[(iii)] $\frac{p_ix_i}{a_i}=\frac{p_jx_j}{a_j}$ for all $i,j \in S^x$,
     \item[(iv)]  $x_i=0$ for all $i \in N \backslash S^x$.
 \end{enumerate}
 \end{definition}

 \noindent A strategy  is semi-balanced if   there  is a subset  of facilities  that gets resources allocated in such a way that the follower becomes indifferent between those facilities (similar to a balanced strategy). Moreover there exists exactly one other facility  that receives  the remaining resources in such a way that its destruction ratio is strictly  less than the destruction ratio of the facilities in the subset. Definition 2 provides the formal definition of a semi-balanced strategy.

\begin{definition}
A strategy $x$ is \emph{semi-balanced} if there exists a set of facilities $S^x \subseteq N$ and a facility $r^x \in N \backslash S^x$  such that
\begin{enumerate}[leftmargin=1.3cm]
\item[(i)]$\sum_{i \in N}{x_i} = R_l$,
    \item[(ii)] $x_i>0$ for all $i \in S^x \cup \{r^x\}$,
    \item[(iii)]  $\frac{p_ix_i}{a_i}=\frac{p_jx_j}{a_j}$ for all $i,j \in S^x$,
    \item[(iv)] $\frac{p_{r^x}x_{r^x}}{a_{r^x}}<\frac{p_ix_i}{a_i}$ for all $i \in S^x$,
    \item[(v)] $x_i=0$ for all $i \in N \backslash (S^x \cup \{r^x\})$.
\end{enumerate}
\end{definition}

\noindent We  illustrate  the definitions of balancedness and semi-balancedness by means of two examples.

\begin{example} \label{example4} Assume the leader applies strategy $x' = (\frac{3}{8},\frac{5}{8},\frac{1}{4},0,3 \frac{3}{4})$.   The resulting destruction ratios $\{\frac{p_ix'_i}{a_i}\}_{i \in N}$ are presented in Table \ref{tableexample4}.
\begin{table}[ht] \centering
\begin{tabular}{lccccc} \hline
Facility $i$  & 1 & 2 & 3 & 4 & 5 \\ \hline
$\frac{p_ix'_i}{a_i} $ & $5$ & $5$ & $5$ &  0 & $5$ \\\hline
\end{tabular}
\caption{The destruction ratios of Example~\ref{example4}}
\label{tableexample4}
\end{table}

\noindent From this table we can conclude that $x'$ is a balanced strategy with destruction ratio $5$ and $S^{x'} = \{1,2,3,5\}$. Another example of a balanced strategy is strategy $x'' =\left(0, \frac{25}{37},\frac{10}{37},0,4\frac{2}{37}\right)$ with destruction ratio $5\frac{15}{37}$ and $S^{x''}=\{2,3,5\}$.
\end{example}

\begin{example} \label{example3} Assume the leader applies strategy $\hat{x} = (\frac{1}{15}, \frac{2}{3},\frac{4}{15},0,4)$.   The resulting destruction ratios $\{\frac{p_i\hat{x}_i}{a_i}\}_{i \in N}$ are presented in Table \ref{tableexample3}. \newpage
\begin{table}[ht] \centering
\begin{tabular}{lccccc} \hline
Facility $i$  & 1 & 2 & 3 & 4 & 5 \\ \hline
$\frac{p_i\hat{x}_i}{a_i} $ & $\frac{8}{9}$ & $5 \frac{1}{3}$ & $5 \frac{1}{3}$ &  0 & $5 \frac{1}{3}$ \\\hline
\end{tabular}
\caption{The destruction ratios of Example~\ref{example3}}
\label{tableexample3}
\end{table}

\noindent From this table we can conclude that $\hat{x}$ is a semi-balanced strategy with destruction ratio $5\frac{1}{3}$, $S^{\hat{x}} = \{2,3,5\}$ and $r^{\hat{x}} = 1$.\end{example}

Now, we  show that  it suffices to consider the class of balanced and semi-balanced strategies as formalized in the following lemma.

\begin{lemma}
\label{towardsgreedycap}
There exists an optimal strategy for the leader that is either balanced or semi-balanced.
\end{lemma}

The following two examples illustrate the main idea of the proof of Lemma~\ref{towardsgreedycap}, namely that each strategy is dominated by an associated balanced or semi-balanced strategy.

\begin{example} \label{example5} Reconsider strategy $x= \left(0,\frac{7}{10},\frac{3}{10},0,4\right)$. Note that from Table~\ref{table2} it follows that this strategy is neither balanced nor semi-balanced. \textcolor{black}{This is also visualized in Figure~\ref{fig:destructionratios}(a).}  From Example \ref{example2} we learned that in this case the follower will fully destroy facility 2 and 3, and he will destroy facility 5 partly.

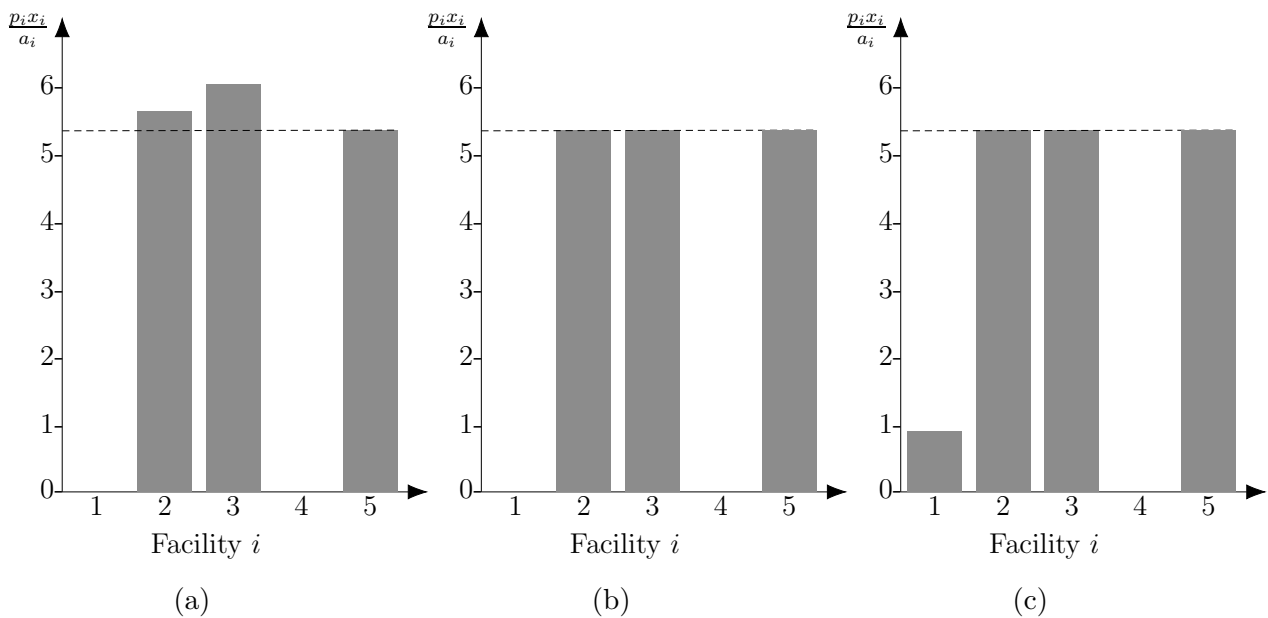
\begin{figure}[ht]
\begin{subfigure}{.32\textwidth}
  \centering
 \scalebox{0.9}{
   \begin{tikzpicture}
   \draw(-0.6,6.8) node {$\frac{p_ix_i}{a_i}$};
    \draw(-0.25,0.05) node {$0$-};
   \draw(-0.25,1) node {$1$-};
   \draw(-0.25,2) node {$2$-};
   \draw(-0.25,3) node {$3$-};
   \draw(-0.25,4) node {$4$-};
      \draw(-0.25,5) node {$5$-};
   \draw(-0.25,6) node {$6$-};
   \draw [-{Latex[length=3mm]}]   (-0.1,0) -- (-0.1,7);
   \draw [-{Latex[length=3mm]}]   (-0.1,0) -- (5.25,0);
   \draw(0.4,-0.2) node {1};
    \draw(1.4,-0.2) node {2};
     \draw(2.4,-0.2) node {3};
      \draw(3.4,-0.2) node {4};
       \draw(4.4,-0.2) node {5};
\path [fill=gray!90!white] (1,0) rectangle (1.8,5.60);
\path [fill=gray!90!white] (2,0) rectangle (2.8,6);
\path [fill=gray!90!white] (4,0) rectangle (4.8,5.33);
 \draw [densely dashed]   (-0.1,5.32) -- (4.8,5.33);
   \draw(2,-0.8) node {Facility $i$};
\end{tikzpicture}}
  \caption{}
  \label{fig:sub-first}
\end{subfigure}
\begin{subfigure}{.32\textwidth}
  \centering
  \scalebox{0.9}{
   \begin{tikzpicture}
   \draw(-0.6,6.8) node {$\frac{p_ix_i}{a_i}$};
   \draw(-0.25,0.05) node {$0$-};
   \draw(-0.25,1) node {$1$-};
   \draw(-0.25,2) node {$2$-};
   \draw(-0.25,3) node {$3$-};
   \draw(-0.25,4) node {$4$-};
      \draw(-0.25,5) node {$5$-};
   \draw(-0.25,6) node {$6$-};
   \draw [-{Latex[length=3mm]}]   (-0.1,0) -- (-0.1,7);
   \draw [-{Latex[length=3mm]}]   (-0.1,0) -- (5.25,0);
   \draw(0.4,-0.2) node {1};
    \draw(1.4,-0.2) node {2};
     \draw(2.4,-0.2) node {3};
      \draw(3.4,-0.2) node {4};
       \draw(4.4,-0.2) node {5};
\path [fill=gray!90!white] (1,0) rectangle (1.8,5.33);
\path [fill=gray!90!white] (2,0) rectangle (2.8,5.33);
 \draw [densely dashed]   (-0.1,5.32) -- (4.8,5.33);
\path [fill=gray!90!white] (4,0) rectangle (4.8,5.33);
   \draw(2,-0.8) node {Facility $i$};
\end{tikzpicture}}
  \caption{}
  \label{fig:sub-second}
\end{subfigure}
\begin{subfigure}{.32\textwidth}
  \centering
  \scalebox{0.9}{
   \begin{tikzpicture}
   \draw(-0.6,6.8) node {$\frac{p_ix_i}{a_i}$};
   \draw(-0.25,0.05) node {$0$-};
   \draw(-0.25,1) node {$1$-};
   \draw(-0.25,2) node {$2$-};
   \draw(-0.25,3) node {$3$-};
   \draw(-0.25,4) node {$4$-};
      \draw(-0.25,5) node {$5$-};
   \draw(-0.25,6) node {$6$-};
   \draw [-{Latex[length=3mm]}]   (-0.1,0) -- (-0.1,7);
   \draw [-{Latex[length=3mm]}]   (-0.1,0) -- (5.25,0);
   \draw(0.4,-0.2) node {1};
    \draw(1.4,-0.2) node {2};
     \draw(2.4,-0.2) node {3};
      \draw(3.4,-0.2) node {4};
       \draw(4.4,-0.2) node {5};
 \path [fill=gray!90!white] (0,0) rectangle (0.8,0.89);
\path [fill=gray!90!white] (1,0) rectangle (1.8,5.33);
\path [fill=gray!90!white] (2,0) rectangle (2.8,5.33);
 \draw [densely dashed]   (-0.1,5.32) -- (4.8,5.33);
\path [fill=gray!90!white] (4,0) rectangle (4.8,5.33);
   \draw(2,-0.8) node {Facility $i$};
\end{tikzpicture}}
  \caption{}
  \label{fig:sub-third}
\end{subfigure}
\caption{Visualization of the destruction ratios of the strategies in Example~\ref{example5}}
\label{fig:destructionratios}
\end{figure}

Suppose now that the leader will adjust his strategy in such a way that the destruction ratios of facility 2 and 3 become equal to the destruction ratio of facility 5, namely $5 \frac{1}{3}$. The leader achieves this by reducing the number of resources allocated to facility 2 and 3. That means, by allocating $\frac{a_2}{p_2} \cdot 5 \tfrac{1}{3} = \frac{2}{3}$ resources to facility 2 and $\frac{a_3}{p_3} \cdot 5 \tfrac{1}{3} = \frac{4}{15}$ resources to facility 3. \textcolor{black}{The destruction ratios of this alternative strategy are visualized in Figure~\ref{fig:destructionratios}(b).} Observe that, for this alternative strategy, the optimal strategy for the follower will not change (since the order $(3,2,5,1,4)$ is still a non-increasing order of the facilities with respect to the destruction ratios). Hence, also the worst case total production after destruction stays the same.

Due to this reduction of resources allocated to facilities 2 and 3, there are some unused resources left. In total there are now  $\frac{7}{10}+\frac{3}{10} - (\frac{2}{3} + \frac{4}{15})=\frac{1}{15}$ unused resources. We can use those resources for increasing the worst case total production after destruction. This can be done by allocating the resources to a facility that will not be destroyed by the follower, so facility 1 or 4. Since facility 1 has the highest production rate, it is most beneficial to allocate these $\frac{1}{15}$ unused resources to facility 1. This results in a destruction ratio of $\frac{p_1}{a_1} \frac{1}{15}=\frac{8}{9}$ for facility 1 and thus the optimal strategy for the follower will not change (since the order $(3,2,5,1,4)$ is still a non-increasing order of the facilities with respect to the destruction ratios). As a result, for this new strategy, which we denote by $\hat{x} = (\frac{1}{15}, \frac{2}{3}, \frac{4}{15},0,4)$, the worst case total production after destruction will increase by $p_1\hat{x}_1$ and thus
$$\mathscr{P}(\hat{x}) = \mathscr{P}(x) + p_1\hat{x}_1 = 1\tfrac{1}{3}+\tfrac{12}{15} = 2 \tfrac{2}{15}.$$
Note that from Example~\ref{example3} it follows that $\hat{x}$ is semi-balanced \textcolor{black}{(which is visualized in Figure~\ref{fig:destructionratios}(c)).}  Hence, we found a semi-balanced strategy $\hat{x}$ that dominates our initial strategy $x$.
\end{example}

The previous example illustrates that a strategy, that is neither balanced nor semi-balanced, is dominated by a semi-balanced strategy. In other cases, it might be that such a strategy is dominated by a balanced strategy. This is illustrated in the following example.

\begin{example}
\label{example5continued}
Reconsider Example~\ref{example5} and suppose that, after the described adjustment procedure to obtain the semi-balanced strategy $\hat{x}$, with $S^{\hat{x}}=\{2,3,5\}$ and $r^{\hat{x}}=\{1\}$, still some unused resources are left (i.e., $R_l>5$). Then, the leader can improve by allocating these unused resources to a facility that will not be destroyed, so facility 1 or 4. It is most beneficial to allocate these unused resources to facility 1. The leader can do so until no unused resources are left or until facility 1 also reaches the destruction ratio of $5 \tfrac{1}{3}$ (cf. Figure~\ref{fig:destructionratios2}(a)). In the latter case, the leader continues with allocating the unused resources to facility 4 until either no unused resources are left (cf. Figure~\ref{fig:destructionratios2}(b)) or until facility 4 also reaches the destruction ratio of $5 \tfrac{1}{3}$. If still some unused resources are left, the leader continues with allocating the resources to all facilities such that the destruction ratios remain the same for every facility (cf. Figure~\ref{fig:destructionratios2}(c)). The leader can do so until no unused resources are left.

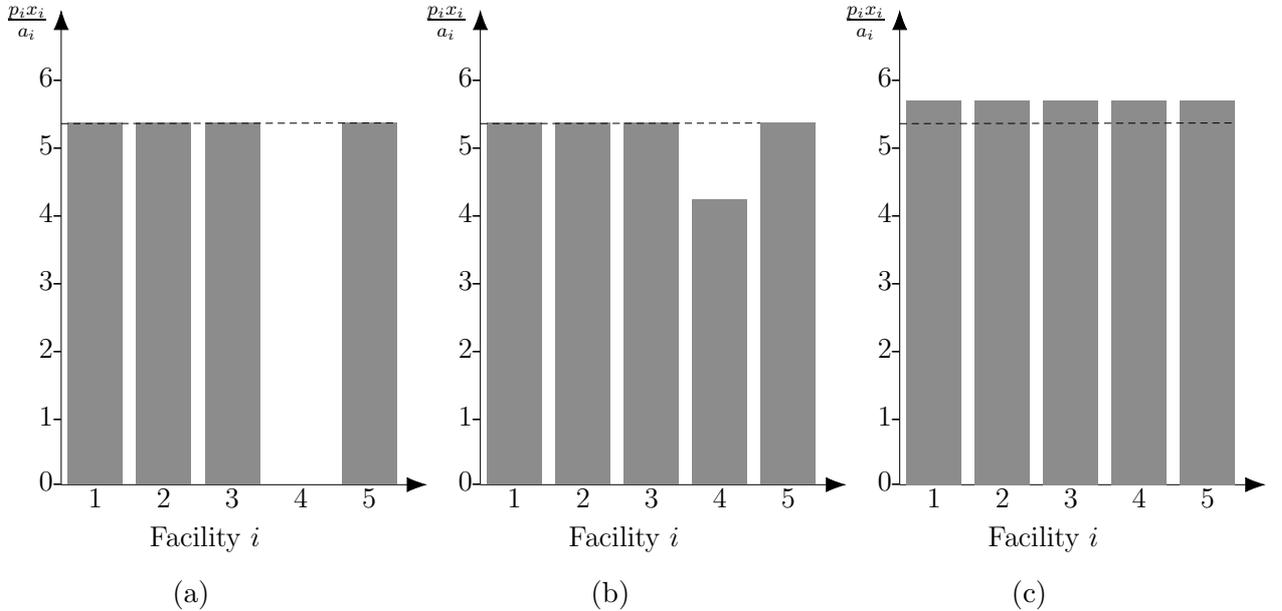
\begin{figure}[ht]
\begin{subfigure}{.32\textwidth}
  \centering
 \scalebox{0.9}{
   \begin{tikzpicture}
   \draw(-0.6,6.8) node {$\frac{p_ix_i}{a_i}$};
    \draw(-0.25,0.05) node {$0$-};
   \draw(-0.25,1) node {$1$-};
   \draw(-0.25,2) node {$2$-};
   \draw(-0.25,3) node {$3$-};
   \draw(-0.25,4) node {$4$-};
      \draw(-0.25,5) node {$5$-};
   \draw(-0.25,6) node {$6$-};
   \draw [-{Latex[length=3mm]}]   (-0.1,0) -- (-0.1,7);
   \draw [-{Latex[length=3mm]}]   (-0.1,0) -- (5.25,0);
   \draw(0.4,-0.2) node {1};
    \draw(1.4,-0.2) node {2};
     \draw(2.4,-0.2) node {3};
      \draw(3.4,-0.2) node {4};
       \draw(4.4,-0.2) node {5};
 \path [fill=gray!90!white] (0,0) rectangle (0.8,5.33);
\path [fill=gray!90!white] (1,0) rectangle (1.8,5.33);
\path [fill=gray!90!white] (2,0) rectangle (2.8,5.33);

\path [fill=gray!90!white] (4,0) rectangle (4.8,5.33);
 \draw [densely dashed]   (-0.1,5.32) -- (4.8,5.33);
   \draw(2,-0.8) node {Facility $i$};
\end{tikzpicture}}
  \caption{}
  \label{fig:sub-first1}
\end{subfigure}
\begin{subfigure}{.32\textwidth}
  \centering
  \scalebox{0.9}{
   \begin{tikzpicture}
   \draw(-0.6,6.8) node {$\frac{p_ix_i}{a_i}$};
   \draw(-0.25,0.05) node {$0$-};
   \draw(-0.25,1) node {$1$-};
   \draw(-0.25,2) node {$2$-};
   \draw(-0.25,3) node {$3$-};
   \draw(-0.25,4) node {$4$-};
      \draw(-0.25,5) node {$5$-};
   \draw(-0.25,6) node {$6$-};
   \draw [-{Latex[length=3mm]}]   (-0.1,0) -- (-0.1,7);
   \draw [-{Latex[length=3mm]}]   (-0.1,0) -- (5.25,0);
   \draw(0.4,-0.2) node {1};
    \draw(1.4,-0.2) node {2};
     \draw(2.4,-0.2) node {3};
      \draw(3.4,-0.2) node {4};
       \draw(4.4,-0.2) node {5};
  \path [fill=gray!90!white] (0,0) rectangle (0.8,5.33);
\path [fill=gray!90!white] (1,0) rectangle (1.8,5.33);
\path [fill=gray!90!white] (2,0) rectangle (2.8,5.33);
\path [fill=gray!90!white] (3,0) rectangle (3.8,4.2);
 \draw [densely dashed]   (-0.1,5.32) -- (4.8,5.33);
\path [fill=gray!90!white] (4,0) rectangle (4.8,5.33);
   \draw(2,-0.8) node {Facility $i$};
\end{tikzpicture}}
  \caption{}
  \label{fig:sub-second1}
\end{subfigure}
\begin{subfigure}{.32\textwidth}
  \centering
  \scalebox{0.9}{
   \begin{tikzpicture}
   \draw(-0.6,6.8) node {$\frac{p_ix_i}{a_i}$};
   \draw(-0.25,0.05) node {$0$-};
   \draw(-0.25,1) node {$1$-};
   \draw(-0.25,2) node {$2$-};
   \draw(-0.25,3) node {$3$-};
   \draw(-0.25,4) node {$4$-};
      \draw(-0.25,5) node {$5$-};
   \draw(-0.25,6) node {$6$-};
   \draw [-{Latex[length=3mm]}]   (-0.1,0) -- (-0.1,7);
   \draw [-{Latex[length=3mm]}]   (-0.1,0) -- (5.25,0);
   \draw(0.4,-0.2) node {1};
    \draw(1.4,-0.2) node {2};
     \draw(2.4,-0.2) node {3};
      \draw(3.4,-0.2) node {4};
       \draw(4.4,-0.2) node {5};
 \path [fill=gray!90!white] (0,0) rectangle (0.8,5.66);
\path [fill=gray!90!white] (1,0) rectangle (1.8,5.66);
\path [fill=gray!90!white] (2,0) rectangle (2.8,5.66);

\path [fill=gray!90!white] (3,0) rectangle (3.8,5.66);
\path [fill=gray!90!white] (4,0) rectangle (4.8,5.66);
 \draw [densely dashed]   (-0.1,5.32) -- (4.8,5.33);
   \draw(2,-0.8) node {Facility $i$};
\end{tikzpicture}}
  \caption{}
  \label{fig:sub-third1}
\end{subfigure}
\caption{Visualization of the destruction ratios of the strategies in Example~\ref{example5continued} }
\label{fig:destructionratios2}
\end{figure}



Observe that for all these three strategies, as illustrated in Figure~\ref{fig:destructionratios2}, the optimal strategy for the follower will not change. As a result, for every possible scenario, we obtain a strategy for the leader that is improved with respect to strategy $x$. Moreover, note that those three strategies are all either balanced (Figure~\ref{fig:destructionratios2}(a) and (c)) or semi-balanced (Figure~\ref{fig:destructionratios2}(b)). Hence, for every possible scenario, we found a (semi-) balanced strategy that dominates our initial strategy $x$.
\end{example}


\subsection{Reducing to the class of balanced strategies}
\label{sectie:bal}

In the previous section, we showed that, for finding an optimal strategy, it suffices to consider the class of balanced and semi-balanced strategies. In this section, we  show that it suffices to consider the class of balanced strategies only. We do so by showing that each semi-balanced strategy is dominated by a balanced strategy. For this, we first consider balanced and semi-balanced strategies in more detail.

The following lemma states the worst case total production after destruction for a balanced strategy. In addition, this lemma also provides the allocation of all resources to the facilities.



\begin{lemma}
\label{allocationbalancedx+balanced}
Let $x$ be a balanced strategy. Then,
\begin{equation*}
\mathscr{P}(x) = \begin{cases}
{\displaystyle \frac{\left ( \sum_{i \in S^x} a_i - R_f \right ) R_l}{\sum_{i \in S^x} \frac{a_i}{p_i}}} & \text{ if } R_f \leq \sum_{i \in S^x}{a_i}, \\
0 & \text{ if } R_f > \sum_{i \in S^x}{a_i},
 \end{cases}
\end{equation*}
and
\begin{equation*}
x_i = \begin{cases}
{\displaystyle \frac{a_iR_l}{p_i \sum_{j \in S^x} \frac{a_j}{p_j}}} & \text{ if } i \in S^x, \\
0 & \text{ if } i \in N \backslash S^x.
 \end{cases}
\end{equation*}
\end{lemma}


From the previous lemma it follows that whenever the leader has access to more resources, ceteris paribus, the worst case total production after destruction will increase with rate $\frac{ \sum_{i \in S^x} a_i - R_f }{\sum_{i \in S^x} \frac{a_i}{p_i}}$. We call this rate the composed net production rate and it represents, for a given balanced strategy, the worst case total production after destruction per invested resource. Note that this rate equals zero in case the follower has enough destructive resources to destroy all facilities of the balanced strategy. The following definition summarizes the above observations.

\begin{definition}
\label{def:avgrate}
The \emph{composed net production rate} is, for all $S \subseteq N$ with $ S \not = \emptyset$, defined by
$$\overline{p}(S)= \max \left \{\frac{\sum_{i \in S} a_i - R_f}{\sum_{i \in S} \frac{a_i}{p_i}},0 \right \},$$
and $\overline{p}(\emptyset)=0$.
\end{definition}

The following corollary follows immediately from Lemma~\ref{allocationbalancedx+balanced} and Definition~\ref{def:avgrate}.
\begin{corollary}    \label{cor:allocationbalancedx+balancedNIEUW}
Let $x$ be a balanced strategy. Then,
    $$\mathscr{P}(x) = \overline{p}(S^x)R_l.$$
\end{corollary}
The following example illustrates the result of the Corollary~\ref{cor:allocationbalancedx+balancedNIEUW}.

\begin{example}
It follows from Corollary~\ref{cor:allocationbalancedx+balancedNIEUW} that the worst case total production after destruction of balanced strategies $x'$ with $S^{x'} = \{1,2,3,5\}$ and $x''$ with $S^{x''} = \{2,3,5\}$ are
\begin{align*}
    \mathscr{P}(x') &= \overline{p}(\{1,2,3,5\})R_l = \max \left \{ \frac{a_1+a_2+a_3+a_5-R_f}{\frac{a_1}{p_1}+\frac{a_2}{p_2}+\frac{a_3}{p_3}+\frac{a_5}{p_5}}, 0 \right \} R_l = 1\tfrac{3}{20} \cdot 5 = 5 \tfrac{3}{4}, \\
    \mathscr{P}(x'') &= \overline{p}(\{2,3,5\})R_l = \max \left \{ \frac{a_2+a_3+a_5-R_f}{\frac{a_2}{p_2}+\frac{a_3}{p_3}+\frac{a_5}{p_5}}, 0 \right \} R_l = \tfrac{10}{37} \cdot 5 = 1 \tfrac{13}{37}. \qedhere
\end{align*}
\end{example}



We now continue with Lemma~\ref{semi-balanced} that states the worst case total production after destruction for a semi-balanced strategy. Note that this lemma also uses the definition of composed net production rate. In addition, this lemma also provides the allocation of all resources to the facilities.



\begin{lemma}
\label{semi-balanced} Let $x$ be a semi-balanced strategy. Then,
\begin{equation*}
\mathscr{P}(x) = \begin{cases}
{\displaystyle \overline{p}(S^x) (R_l-x_{r^{x}}) + p_{r^{x}} x_{r^{x}}} & \text{ if } R_f \leq  \sum_{i \in S^x}{a_i} , \\
{\displaystyle \frac{p_{r^{x}}x_{r^{x}}}{a_{r^{x}}} \left( \sum_{i \in S^x \cup \{r^{x}\}} a_i-R_f\right)} & \text{ if } \sum_{i \in S^x} a_i <  R_f < \sum_{i \in S^x \cup \{r^{x}\}} a_i, \\
0 & \text{ if } R_f \geq \sum_{i \in S^x \cup \{r^{x}\}} a_i,
 \end{cases}
\end{equation*}
and
\begin{equation*}
x_i = \begin{cases}
{\displaystyle \frac{a_i(R_l-x_{r^{x}})}{p_i \sum_{j \in S^x} \frac{a_j}{p_j}}} & \text{ if } i \in S^x, \\
0 & \text{ if } i \in N \backslash (S^x \cup r^{x}).
 \end{cases}
\end{equation*}
\end{lemma}


The following example illustrates  the previous lemma.




\begin{example}
Reconsider the semi-balanced strategy  $\hat{x}$ with $S^{\hat{x}}=\{2,3,5\}$ and $r^{\hat{x}}=\{1\}$. We have $R_f=1 \frac{3}{4} \leq  1+ \frac{1}{4}  + \frac{3}{4}= a_2+a_3+a_5 = \sum_{i \in S^{\hat{x}}}{a_i}$ and thus, using Lemma~\ref{semi-balanced}, the worst case total production after destruction of this semi-balanced strategy equals
\begin{align*}
\mathscr{P}(\hat{x}) &= \overline{p}(\{2,3,5\}) (R_l-\hat{x}_1) + p_{1} \hat{x}_{1} \\
&= \max \left \{ \frac{a_2+a_3+a_5-R_f}{\frac{a_2}{p_2}+\frac{a_3}{p_3}+\frac{a_5}{p_5}}, 0 \right \} (R_l-\hat{x}_1) + p_{1} \hat{x}_{1} \\
&=\tfrac{10}{37} (5- \tfrac{1}{15}) + 12 \cdot \tfrac{1}{15} \\
&=2 \tfrac{2}{15}.
\end{align*}
Note  we  already calculated this worst case total production after destruction in Example~\ref{example5}.
\end{example}



By using Corollary~\ref{cor:allocationbalancedx+balancedNIEUW} and Lemma~\ref{semi-balanced}, we are now able to show that it suffices to consider the class of balanced strategies as formalized in the following lemma.

\begin{lemma}
\label{t:nofractions}
There exists an optimal strategy for the leader that is balanced.
\end{lemma}

The main idea of the proof of Lemma~\ref{t:nofractions} is that each semi-balanced strategy is dominated by a balanced strategy. Namely, a semi-balanced strategy $x$ is either dominated by the balanced strategy $S^x$ or the balanced strategy $S^x \cup \{r^x\}$. The following example illustrates why this domination always holds. In particular, we  illustrate in this example how the proof makes use of the composed net production rate and the results of Corollary~\ref{cor:allocationbalancedx+balancedNIEUW} and Lemma~\ref{semi-balanced}.

\begin{example} \label{examplezoveel}
Reconsider the semi-balanced strategy  $\hat{x}$ with $S^{\hat{x}}=\{2,3,5\}$ and $r^{\hat{x}}=\{1\}$. Observe that the production rate of  $r^{\hat{x}}$ is higher than the composed net production rate of $S^{\hat{x}}$, i.e.,
\begin{equation}
\label{eq:p1versuspstreep}
    p_1=12 > \tfrac{10}{37} = \overline{p}(\{2,3,5\}).
\end{equation}
From Lemma \ref{semi-balanced} we can conclude that exchanging resources from $S^{\hat{x}}$ to facility 1, as long as the follower's optimal strategy stays the same, is beneficial. In  other words, in order to improve the strategy for the leader, one can increase $\hat{x}_1$ (and thus decrease $\hat{x}_2$, $\hat{x}_3$ and $\hat{x}_5$) such that the destruction ratio of facility 1 becomes equal to the destruction ratios of  $S^{\hat{x}}$. This results in the balanced strategy $x'$ with $S^{x'} = \{1,2,3,5\}$, which dominates the semi-balanced strategy $\hat{x}$.

Note that in case the inequality in~\eqref{eq:p1versuspstreep} would be reversed, then exchanging resources from facility 1 to $S^{\hat{x}}$ would be beneficial. As a result, the balanced strategy $x''$ with $S^{x''} = \{2,3,5\}$ would dominate the semi-balanced strategy $\hat{x}$.
\end{example}

\subsection{Reducing to the class of seried-balanced strategies}
\label{sectie:seried}

In this section, we  show that, for finding an optimal strategy, it suffices to consider a proper subset  of the class of balanced strategies: the class of seried-balanced strategies. A balanced strategy is called seried-balanced if only facilities up to a certain threshold index get resources allocated. Definition~\ref{def:seried-balanced} provides the formal definition of a seried-balanced strategy.




\begin{definition}\label{def:seried-balanced}
A balanced strategy $x \in \mathscr{X}$ is \emph{seried-balanced} if
$$S^{x} = \{1,2,\ldots,s^x\},$$
where $s^x$ the highest indexed facility in $S^x$, i.e., $s^x = \max \{i~\vert~i \in S^x\}$.
\end{definition}
Recall that the production rates are in non-increasing order (i.e., $p_1 \geq p_2 \geq \hdots \geq p_n$). As a consequence, a seried-balanced strategy always includes the facilities with the highest production rates. The following example illustrates the definition of a seried-balanced strategy.

\begin{example}
\label{example:notseried}
The set of all seried-balanced strategies is given by
\begin{equation*}
\left \{x \in \mathscr{X}~|~\textrm{ $x$ is balanced, }S^x \in \{\{1\}, \{1,2\}, \{1,2,3\}, \{1,2,3,4\}, \{1,2,3,4,5\}\} \right \}.
\end{equation*}
Hence, balanced strategy $x'$ with $S^{x'} = \{1,2,3,5\}$  is not seried-balanced, since facility $4$ does not get any resources allocated while facility $5$ does.
\end{example}

In order to prove that we can restrict our attention to the class of seried-balanced strategies, we show that every balanced strategy is dominated by a seried-balanced strategy. We show this domination by adding facilities to a balanced strategy and/or removing facilities from a balanced strategy until the strategy becomes seried-balanced. The following lemma states when it is beneficial (or not) to add a facility to (or remove a facility from) a balanced strategy.


\begin{lemma}
\label{veiligtoevoegen}
Let $x$ and $x'$  be two balanced strategies such that $S^x \cap S^{x'} = \emptyset$. Moreover, let $x''$ be the balanced strategy such that $S^{x''}=S^x \cup S^{x'}$.
\begin{enumerate}[label=(\roman*)]
\item If $p_i \geq \overline{p}(S^x)$ for all $i \in S^{x'}$, then $\overline{p}(S^{x''}) \geq \overline{p}(S^x)$;
\item If $p_i \leq \overline{p}(S^x)$ for all $i \in S^{x'}$, then $\overline{p}(S^{x''}) \leq \overline{p}(S^x)$;
\item If $p_i \geq \overline{p}(S^{x''})$ for all $i \in S^{x'}$, then $\overline{p}(S^{x''}) \geq \overline{p}(S^x)$;
\item If $p_i \leq \overline{p}(S^{x''})$ for all $i \in S^{x'}$, then $\overline{p}(S^{x''}) \leq \overline{p}(S^x)$.
\end{enumerate}
\end{lemma}

By exploiting Lemma \ref{veiligtoevoegen}, we are able to show that it suffices to consider the class of seried-balanced strategies. This is formalized in the following lemma.

\begin{lemma}
\label{bijnaklaarstelling}
There exists an optimal strategy for the leader that is seried-balanced.
\end{lemma}

The following example illustrates the main idea of the proof of Lemma~\ref{bijnaklaarstelling}. Namely, how one should add facilities to a balanced strategy and/or remove facilities from a balanced strategy until a seried-balanced strategy is obtained that dominates the initial balanced strategy.


\begin{example} \label{examplezoveel2}
Reconsider  the balanced strategy $x'$ with $S^{x'}=\{1,2,3,5\}$  and recall that this strategy is not seried-balanced. Note that the production rate of facility 4 is higher than the composed net production rate of $S^{x'}$, i.e.,
$$
    p_4=2 > 1 \tfrac{3}{20} = \overline{p}(\{1,2,3,5\}).
$$
Therefore, it follows from Lemma~\ref{veiligtoevoegen}\emph{(i)} that $\overline{p}(\{1,2,3,4,5\}) \geq \overline{p}(\{1,2,3,5\})$. As a result, due to Corollary~\ref{cor:allocationbalancedx+balancedNIEUW}, the strategy that is balanced over all facilities dominates balanced strategy $x'$.

Next, note that the production rate of facility 5 is lower than the composed net production rate of all facilities together, i.e.,
$$p_5=1 < 1 \tfrac{13}{30} = \overline{p}(\{1,2,3,4,5\}).
$$
Therefore, it follows from Lemma~\ref{veiligtoevoegen}\emph{(iv)} that $\overline{p}(\{1,2,3,4\}) \geq \overline{p}(\{1,2,3,4,5\})$. As a result, due to Corollary~\ref{cor:allocationbalancedx+balancedNIEUW}, the seried-balanced strategy $x^*$  with $S^{x^*}=\{1,2,3,4\}$ dominates the strategy that is balanced over all facilities.

In summary, the seried-balanced strategy $x^*$ dominates the balanced strategy $x'$, which concludes this example.
 \end{example}

\subsection{Algorithm for finding an optimal strategy for the leader}
\label{sectie:algorithm}
In this section, we provide a linear time algorithm that finds an optimal strategy for the leader. For this, we first make the following three  observations.
\begin{itemize}
    \item For finding an optimal strategy, it suffices to consider the class of seried-balanced strategies (Section~\ref{sectie:seried}).
    \item For a Stackelberg production game with $n$ facilities, there are $n$ seried-balanced strategies (Definition~\ref{def:seried-balanced} and Example \ref{example:notseried}).
    \item For each seried-balanced strategy, the worst case total production after destruction equals the composed net production rate times the leader's amount of resources (Corollary~\ref{cor:allocationbalancedx+balancedNIEUW}).
    \end{itemize}
From these observations it follows that selecting, from all the $n$ seried-balanced strategies, a strategy with the highest composed net production rate results in an optimal strategy. One can do so by calculating the composed net production rate of all these strategies.



We now explain why there is no need to calculate the composed net production rate of all the $n$ seried-balanced strategies. First, recall that the production rates are in non-increasing order (i.e., $p_1 \geq p_2 \geq \ldots \geq p_n)$. Moreover, from Lemma \ref{veiligtoevoegen}($i$) and ($ii$), we learned that adding a facility to a balanced strategy is beneficial (not beneficial), if the production rate of that facility is higher (lower) than  the composed net production rate of that balanced strategy. Hence, if we start with calculating the composed net production rate of facility 1, then calculate the composed net production rate of facility 1 and 2, and so on,
we can stop once the production rate of a facility is at most the composed net production rate of all its preceding facilities. Then, the seried-balanced strategy that corresponds to the last considered composed net production rate is an optimal strategy. In summary, there is no need to calculate the composed net production rate of all the $n$  seried-balanced strategies.

Using the conclusions from above, we present in Algorithm~\ref{alg:algorithm} the pseudocode of a procedure that identifies an optimal strategy for the leader.

\begin{center}
\begin{minipage}{0.85\linewidth}
\begin{algorithm}[H]
\label{alg:algorithm}
\SetAlgoLined
\SetKwInput{KwOutput}{Output}
\SetKwInput{KwInput}{Input}
 \KwInput{Stackelberg production game $\left ( N, \{a_i\}_{i \in N}, \{p_i\}_{i \in N}, R_l, R_f \right ) $}
\KwOutput{an optimal strategy $x^*$}  \smallskip
$S:= \{1\} $; \\
$i:=2$; \\

\While{$p_i > \max \left \{\frac{\sum_{i \in S} a_i - R_f}{\sum_{i \in S} \frac{a_i}{p_i}},0 \right \}$\emph{\textbf{ and }}$i \leq n$\smallskip}{ \smallskip

$S:=S \cup \{i\}$;  \medskip


$i:=i + 1 $\; }
\For{$j=1$\emph{\textbf{ to }}$i-1$}{$x^*_j:={\displaystyle \frac{a_iR_l}{p_i \sum_{j \in S} \frac{a_j}{p_j}}}$;}
\For{$j=i$\emph{\textbf{ to }}$n$}{$x^*_j:=0$;}
 \caption{Finding an optimal strategy}
\end{algorithm}
\end{minipage} \end{center}

 The following theorem states that Algorithm~1 always finds an optimal strategy.


\begin{theorem}
\label{algorithm} Algorithm~\ref{alg:algorithm} is a linear time algorithm for finding an optimal strategy for the leader. \end{theorem}


We conclude this section with an example that explains Algorithm~\ref{alg:algorithm} step by step.




 \begin{example}
\label{example:algorithm}
In the initialization step we set $S:=\{1\}$. We move on with the following iterative steps.

\begin{itemize}
\item \emph{Iteration $i=2$}: Since $p_2 = 8 > 0 = \max \left \{\frac{\sum_{i=1}^{2}{a_i}- R_f}{ \sum_{i=1}^2 \frac{a_i}{p_i}},0\right \}$, it follows from Lemma~\ref{veiligtoevoegen}($i$) that it is beneficial to add facility 2 to the seried-balanced strategy and thus we set $S:=S \cup \{2\} = \{1,2\}$.
\item \emph{Iteration $i=3$}: Since $p_3 = 5 > \tfrac{3}{4} = \max \left \{\frac{\sum_{i=1}^{3}{a_i}- R_f}{ \sum_{i=1}^3 \frac{a_i}{p_i}},0\right \}$, we set $S:=S \cup \{3\} = \{1,2,3\}$.
\item \emph{Iteration $i=4$}: Since $p_4 = 2 >1\tfrac{3}{5} =  \max \left \{\frac{\sum_{i=1}^4 a_i - R_f}{\sum_{i=1}^4\frac{a_i}{p_i}},0 \right \}$, we set $S:=S \cup \{4\} = \{1,2,3,4\}$.
\item \emph{Iteration $i=5$}: Since $p_5 = 1 \not > 1\tfrac{13}{15}
= \max \left \{\frac{\sum_{i=1}^{5}{a_i}- R_f}{ \sum_{i=1}^5 \frac{a_i}{p_i}},0\right \}$, it follows from Lemma~\ref{veiligtoevoegen}($ii$) that it is not beneficial to add facility 5 to the seried-balanced strategy and thus this iteration will not be considered and $S$ will not be updated.
\end{itemize}
Hence, the algorithm returns the seried-balanced strategy $x^* = (\tfrac{1}{2},\tfrac{5}{6},\tfrac{1}{3},3\tfrac{1}{3},0)$. By Theorem~\ref{algorithm} it follows that this strategy is an optimal strategy for the leader.
Note that under this optimal strategy $x^*$, we have $\mathscr{P}(x^*) = 1 \tfrac{13}{15} \cdot 5 = 9 \tfrac{1}{3}.$
\end{example}
\newpage

\section{Conclusions} \label{sectie: conclusions}

Inspired by a military context, we studied a Stackelberg production game, where a leader wants to maximize his production and a  follower tries to destroy this production as much as possible. This paper focused on identifying the leader's and follower's optimal strategies. In particular, we showed that a leader's optimal production strategy can be found in the class of seried-balanced
strategies. Moreover, we presented a linear time algorithm that finds an optimal production strategy in this class.

For future research, we identify three research directions. These directions are based on our modelling assumptions. Firstly, we implicitly assumed in our model that each production facility has ample production capacity (i.e., the leader can, for example, allocate all his resources to one production facility). In practice, however, production facilities are sometimes restricted by their production capacities. This calls for studying our model with a maximum production capacity constraint per production facility. By incorporating this constraint, balanced strategies are not always feasible anymore and thus it is unclear whether each production strategy is still dominated by a balanced strategy.


Secondly, we assumed in our model that the difficulty of destroying a production facility is independent of the number of military assets produced per production facility. However, in practice, there might be a relationship between these two components. Including such a dependency into our model forms a second direction for future research. For this extension, the leader's set of production strategies will depend on the decision of the follower. \citet{fischetti2019interdiction} indicate that incorporating such a dependency might make a Stackelberg game computationally hard. It would be interesting to see whether this is also the case for our Stackelberg production game.




Finally, we assumed in our model that the number of military assets produced is linear in the amount of resources allocated. However, in practice, production quantities might be affected by (dis)economies of scale.
As a consequence, we suggest to study an extension of our model, where the number of military assets produced is non-linear in the amount of resources invested. Since the proof of Lemma \ref{bijnaklaarstelling} uses this linearity, it is unlikely that, for such an extension, we can still restrict our attention to the class of seried-balanced strategies.


\subsection*{Author contribution statement}

L. Schlicher is the corresponding author and initiator of this project. In the early research process, H. Blok and L. Schlicher developed the model and derived the corresponding results. Later on, M. Musegaas also became involved in the project. Together with L. Schlicher, they contributed equally to the writing process.

\bibliographystyle{apacite}
\bibliography{reference}

\begin{thebibliography}{}

\bibitem [\protect \citeauthoryear {%
Bier%
, Oliveros%
\BCBL {}\ \BBA {} Samuelson%
}{%
Bier%
\ \protect \BOthers {.}}{%
{\protect \APACyear {2007}}%
}]{%
bier2007choosing}
\APACinsertmetastar {%
bier2007choosing}%
\begin{APACrefauthors}%
Bier, V.%
, Oliveros, S.%
\BCBL {}\ \BBA {} Samuelson, L.%
\end{APACrefauthors}%
\unskip\
\newblock
\APACrefYearMonthDay{2007}{}{}.
\newblock
{\BBOQ}\APACrefatitle {Choosing what to protect: Strategic defensive allocation
  against an unknown attacker} {Choosing what to protect: Strategic defensive
  allocation against an unknown attacker}.{\BBCQ}
\newblock
\APACjournalVolNumPages{Journal of Public Economic Theory}{9}{4}{563--587}.
\PrintBackRefs{\CurrentBib}

\bibitem [\protect \citeauthoryear {%
Cappanera%
\ \BBA {} Scaparra%
}{%
Cappanera%
\ \BBA {} Scaparra%
}{%
{\protect \APACyear {2011}}%
}]{%
cappanera2011optimal}
\APACinsertmetastar {%
cappanera2011optimal}%
\begin{APACrefauthors}%
Cappanera, P.%
\BCBT {}\ \BBA {} Scaparra, M.%
\end{APACrefauthors}%
\unskip\
\newblock
\APACrefYearMonthDay{2011}{}{}.
\newblock
{\BBOQ}\APACrefatitle {Optimal allocation of protective resources in
  shortest-path networks} {Optimal allocation of protective resources in
  shortest-path networks}.{\BBCQ}
\newblock
\APACjournalVolNumPages{Transportation Science}{45}{1}{64--80}.
\PrintBackRefs{\CurrentBib}

\bibitem [\protect \citeauthoryear {%
Dempe%
}{%
Dempe%
}{%
{\protect \APACyear {2002}}%
}]{%
dempe2002}
\APACinsertmetastar {%
dempe2002}%
\begin{APACrefauthors}%
Dempe, S.%
\end{APACrefauthors}%
\unskip\
\newblock
\APACrefYear{2002}.
\newblock
\APACrefbtitle {Foundations of bilevel programming} {Foundations of bilevel
  programming}.
\newblock
\APACaddressPublisher{}{Springer Science \& Business Media}.
\PrintBackRefs{\CurrentBib}

\bibitem [\protect \citeauthoryear {%
Fischetti%
, Ljubi{\'c}%
, Monaci%
\BCBL {}\ \BBA {} Sinnl%
}{%
Fischetti%
\ \protect \BOthers {.}}{%
{\protect \APACyear {2019}}%
}]{%
fischetti2019interdiction}
\APACinsertmetastar {%
fischetti2019interdiction}%
\begin{APACrefauthors}%
Fischetti, M.%
, Ljubi{\'c}, I.%
, Monaci, M.%
\BCBL {}\ \BBA {} Sinnl, M.%
\end{APACrefauthors}%
\unskip\
\newblock
\APACrefYearMonthDay{2019}{}{}.
\newblock
{\BBOQ}\APACrefatitle {Interdiction games and monotonicity, with application to
  knapsack problems} {Interdiction games and monotonicity, with application to
  knapsack problems}.{\BBCQ}
\newblock
\APACjournalVolNumPages{INFORMS Journal on Computing}{31}{2}{390--410}.
\PrintBackRefs{\CurrentBib}

\bibitem [\protect \citeauthoryear {%
Gutin%
, Kuhn%
\BCBL {}\ \BBA {} Wiesemann%
}{%
Gutin%
\ \protect \BOthers {.}}{%
{\protect \APACyear {2014}}%
}]{%
gutin2014interdiction}
\APACinsertmetastar {%
gutin2014interdiction}%
\begin{APACrefauthors}%
Gutin, E.%
, Kuhn, D.%
\BCBL {}\ \BBA {} Wiesemann, W.%
\end{APACrefauthors}%
\unskip\
\newblock
\APACrefYearMonthDay{2014}{}{}.
\newblock
{\BBOQ}\APACrefatitle {Interdiction games on Markovian PERT networks}
  {Interdiction games on markovian pert networks}.{\BBCQ}
\newblock
\APACjournalVolNumPages{Management Science}{61}{5}{999--1017}.
\PrintBackRefs{\CurrentBib}

\bibitem [\protect \citeauthoryear {%
Hausken%
\ \BBA {} Zhuang%
}{%
Hausken%
\ \BBA {} Zhuang%
}{%
{\protect \APACyear {2011}}%
}]{%
hausken2011defending}
\APACinsertmetastar {%
hausken2011defending}%
\begin{APACrefauthors}%
Hausken, K.%
\BCBT {}\ \BBA {} Zhuang, J.%
\end{APACrefauthors}%
\unskip\
\newblock
\APACrefYearMonthDay{2011}{}{}.
\newblock
{\BBOQ}\APACrefatitle {Defending against a terrorist who accumulates resources}
  {Defending against a terrorist who accumulates resources}.{\BBCQ}
\newblock
\APACjournalVolNumPages{Military Operations Research}{}{}{21--39}.
\PrintBackRefs{\CurrentBib}

\bibitem [\protect \citeauthoryear {%
Israeli%
\ \BBA {} Wood%
}{%
Israeli%
\ \BBA {} Wood%
}{%
{\protect \APACyear {2002}}%
}]{%
israeli2002shortest}
\APACinsertmetastar {%
israeli2002shortest}%
\begin{APACrefauthors}%
Israeli, E.%
\BCBT {}\ \BBA {} Wood, R.%
\end{APACrefauthors}%
\unskip\
\newblock
\APACrefYearMonthDay{2002}{}{}.
\newblock
{\BBOQ}\APACrefatitle {Shortest-path network interdiction} {Shortest-path
  network interdiction}.{\BBCQ}
\newblock
\APACjournalVolNumPages{Networks: An International Journal}{40}{2}{97--111}.
\PrintBackRefs{\CurrentBib}

\bibitem [\protect \citeauthoryear {%
Jiang%
\ \BBA {} Liu%
}{%
Jiang%
\ \BBA {} Liu%
}{%
{\protect \APACyear {2018}}%
}]{%
jiang2018multi}
\APACinsertmetastar {%
jiang2018multi}%
\begin{APACrefauthors}%
Jiang, J.%
\BCBT {}\ \BBA {} Liu, X.%
\end{APACrefauthors}%
\unskip\
\newblock
\APACrefYearMonthDay{2018}{}{}.
\newblock
{\BBOQ}\APACrefatitle {Multi-objective Stackelberg game model for water supply
  networks against interdictions with incomplete information} {Multi-objective
  stackelberg game model for water supply networks against interdictions with
  incomplete information}.{\BBCQ}
\newblock
\APACjournalVolNumPages{European Journal of Operational
  Research}{266}{3}{920--933}.
\PrintBackRefs{\CurrentBib}

\bibitem [\protect \citeauthoryear {%
Liberatore%
, Scaparra%
\BCBL {}\ \BBA {} Daskin%
}{%
Liberatore%
\ \protect \BOthers {.}}{%
{\protect \APACyear {2011}}%
}]{%
liberatore2011analysis}
\APACinsertmetastar {%
liberatore2011analysis}%
\begin{APACrefauthors}%
Liberatore, F.%
, Scaparra, M.%
\BCBL {}\ \BBA {} Daskin, M.%
\end{APACrefauthors}%
\unskip\
\newblock
\APACrefYearMonthDay{2011}{}{}.
\newblock
{\BBOQ}\APACrefatitle {Analysis of facility protection strategies against an
  uncertain number of attacks: The stochastic R-interdiction median problem
  with fortification} {Analysis of facility protection strategies against an
  uncertain number of attacks: The stochastic r-interdiction median problem
  with fortification}.{\BBCQ}
\newblock
\APACjournalVolNumPages{Computers and Operations Research}{38}{1}{357--366}.
\PrintBackRefs{\CurrentBib}

\bibitem [\protect \citeauthoryear {%
Morton%
, Pan%
\BCBL {}\ \BBA {} Saeger%
}{%
Morton%
\ \protect \BOthers {.}}{%
{\protect \APACyear {2007}}%
}]{%
morton2007models}
\APACinsertmetastar {%
morton2007models}%
\begin{APACrefauthors}%
Morton, D\BPBI P.%
, Pan, F.%
\BCBL {}\ \BBA {} Saeger, K.%
\end{APACrefauthors}%
\unskip\
\newblock
\APACrefYearMonthDay{2007}{}{}.
\newblock
{\BBOQ}\APACrefatitle {Models for nuclear smuggling interdiction} {Models for
  nuclear smuggling interdiction}.{\BBCQ}
\newblock
\APACjournalVolNumPages{IIE Transactions}{39}{1}{3--14}.
\PrintBackRefs{\CurrentBib}

\bibitem [\protect \citeauthoryear {%
Powell%
}{%
Powell%
}{%
{\protect \APACyear {2009}}%
}]{%
powell2009sequential}
\APACinsertmetastar {%
powell2009sequential}%
\begin{APACrefauthors}%
Powell, R.%
\end{APACrefauthors}%
\unskip\
\newblock
\APACrefYearMonthDay{2009}{}{}.
\newblock
{\BBOQ}\APACrefatitle {Sequential, nonzero-sum “Blotto”: Allocating
  defensive resources prior to attack} {Sequential, nonzero-sum “blotto”:
  Allocating defensive resources prior to attack}.{\BBCQ}
\newblock
\APACjournalVolNumPages{Games and Economic Behavior}{67}{2}{611--615}.
\PrintBackRefs{\CurrentBib}

\bibitem [\protect \citeauthoryear {%
Scaparra%
\ \BBA {} Church%
}{%
Scaparra%
\ \BBA {} Church%
}{%
{\protect \APACyear {2008}}%
{\protect \APACexlab {{\protect \BCnt {1}}}}}]{%
scaparra2008bilevel}
\APACinsertmetastar {%
scaparra2008bilevel}%
\begin{APACrefauthors}%
Scaparra, M.%
\BCBT {}\ \BBA {} Church, R.%
\end{APACrefauthors}%
\unskip\
\newblock
\APACrefYearMonthDay{2008{\protect \BCnt {1}}}{}{}.
\newblock
{\BBOQ}\APACrefatitle {A bilevel mixed-integer program for critical
  infrastructure protection planning} {A bilevel mixed-integer program for
  critical infrastructure protection planning}.{\BBCQ}
\newblock
\APACjournalVolNumPages{Computers and Operations Research}{35}{6}{1905--1923}.
\PrintBackRefs{\CurrentBib}

\bibitem [\protect \citeauthoryear {%
Scaparra%
\ \BBA {} Church%
}{%
Scaparra%
\ \BBA {} Church%
}{%
{\protect \APACyear {2008}}%
{\protect \APACexlab {{\protect \BCnt {2}}}}}]{%
scaparra2008exact}
\APACinsertmetastar {%
scaparra2008exact}%
\begin{APACrefauthors}%
Scaparra, M.%
\BCBT {}\ \BBA {} Church, R.%
\end{APACrefauthors}%
\unskip\
\newblock
\APACrefYearMonthDay{2008{\protect \BCnt {2}}}{}{}.
\newblock
{\BBOQ}\APACrefatitle {An exact solution approach for the interdiction median
  problem with fortification} {An exact solution approach for the interdiction
  median problem with fortification}.{\BBCQ}
\newblock
\APACjournalVolNumPages{European Journal of Operational
  Research}{189}{1}{76--92}.
\PrintBackRefs{\CurrentBib}

\bibitem [\protect \citeauthoryear {%
Shan%
\ \BBA {} Zhuang%
}{%
Shan%
\ \BBA {} Zhuang%
}{%
{\protect \APACyear {2013}}%
}]{%
shan2013hybrid}
\APACinsertmetastar {%
shan2013hybrid}%
\begin{APACrefauthors}%
Shan, X.%
\BCBT {}\ \BBA {} Zhuang, J.%
\end{APACrefauthors}%
\unskip\
\newblock
\APACrefYearMonthDay{2013}{}{}.
\newblock
{\BBOQ}\APACrefatitle {Hybrid defensive resource allocations in the face of
  partially strategic attackers in a sequential defender--attacker game}
  {Hybrid defensive resource allocations in the face of partially strategic
  attackers in a sequential defender--attacker game}.{\BBCQ}
\newblock
\APACjournalVolNumPages{European Journal of Operational
  Research}{228}{1}{262--272}.
\PrintBackRefs{\CurrentBib}

\bibitem [\protect \citeauthoryear {%
Starita%
\ \BBA {} Scaparra%
}{%
Starita%
\ \BBA {} Scaparra%
}{%
{\protect \APACyear {2016}}%
}]{%
starita2016optimizing}
\APACinsertmetastar {%
starita2016optimizing}%
\begin{APACrefauthors}%
Starita, S.%
\BCBT {}\ \BBA {} Scaparra, M.%
\end{APACrefauthors}%
\unskip\
\newblock
\APACrefYearMonthDay{2016}{}{}.
\newblock
{\BBOQ}\APACrefatitle {Optimizing dynamic investment decisions for railway
  systems protection} {Optimizing dynamic investment decisions for railway
  systems protection}.{\BBCQ}
\newblock
\APACjournalVolNumPages{European Journal of Operational
  Research}{248}{2}{543--557}.
\PrintBackRefs{\CurrentBib}

\bibitem [\protect \citeauthoryear {%
Washburn%
\ \BBA {} Wood%
}{%
Washburn%
\ \BBA {} Wood%
}{%
{\protect \APACyear {1995}}%
}]{%
washburn1995two}
\APACinsertmetastar {%
washburn1995two}%
\begin{APACrefauthors}%
Washburn, A.%
\BCBT {}\ \BBA {} Wood, K.%
\end{APACrefauthors}%
\unskip\
\newblock
\APACrefYearMonthDay{1995}{}{}.
\newblock
{\BBOQ}\APACrefatitle {Two-person zero-sum games for network interdiction}
  {Two-person zero-sum games for network interdiction}.{\BBCQ}
\newblock
\APACjournalVolNumPages{Operations research}{43}{2}{243--251}.
\PrintBackRefs{\CurrentBib}

\bibitem [\protect \citeauthoryear {%
Zhuang%
\ \BBA {} Bier%
}{%
Zhuang%
\ \BBA {} Bier%
}{%
{\protect \APACyear {2011}}%
}]{%
zhuang2011secrecy}
\APACinsertmetastar {%
zhuang2011secrecy}%
\begin{APACrefauthors}%
Zhuang, J.%
\BCBT {}\ \BBA {} Bier, V.%
\end{APACrefauthors}%
\unskip\
\newblock
\APACrefYearMonthDay{2011}{}{}.
\newblock
{\BBOQ}\APACrefatitle {Secrecy and deception at equilibrium, with applications
  to anti-terrorism resource allocation} {Secrecy and deception at equilibrium,
  with applications to anti-terrorism resource allocation}.{\BBCQ}
\newblock
\APACjournalVolNumPages{Defence and Peace Economics}{22}{1}{43--61}.
\PrintBackRefs{\CurrentBib}

\end{thebibliography}

 \newpage

\section*{Appendix}

In this appendix, we provide the proofs of all lemmas and theorems stated in this paper. For notional convenience in the proofs of Lemma~\ref{towardsgreedycap} and~\ref{semi-balanced}, we rewrite~\eqref{vgl:Loe} as
\begin{equation}
\label{eq:defPxyomega}
\mathscr{P}(x) = x_{q^{\omega}}p_{q^{\omega}} \left (1-\frac{y^{\omega}_{q^{\omega}}}{a_{q^{\omega}}} \right ) + \sum_{i \in N \backslash A^{\omega}}{x_ip_i}.
\end{equation}
Note that this equivalence follows directly from the definition of $A^{\omega}$ and $q^{\omega}$. \bigskip

\noindent \underline{\textbf{Proof of Theorem \ref{thm:optimalstrategyfollower}}}

 \begin{proof}
 The aim of the follower is to minimize the leader's total production after destruction, i.e.,
 \begin{align*} \min_{y \in \mathscr{Y}}\left\{ \mathscr{P}(x,y)\right\}
&=\min_{y \in \mathscr{Y}} \left \{\sum_{i \in N}{x_ip_i \left (1- \frac{y_i}{a_i}  \right )} \right \} \\
 &= \sum_{i \in N} x_i p_i - \max_{y \in \mathscr{Y}} \left\{ \sum_{i \in N} \frac{p_ix_i}{a_i} y_i\right\}. \end{align*}

 \noindent Note that the second equality holds because $\sum_{i \in N} x_ip_i$ is a constant and thus it can be taken outside the minimum. The second equality also follows from the fact that minimizing a function over its argument is equivalent to maximizing that function over the same argument with a sign change.

 Now, let $\omega \in \Omega^x$. Observe that, since the ratio $\frac{p_ix_i}{a_i}$ is a constant for all $i \in N$ in the maximum, it is optimal for the follower to consider the facilities in non-increasing order with respect to ratio $\frac{p_ix_i}{a_i}$, i.e., to consider the facilities in order $\omega$.
Hence, the follower allocates his destructive resources facility by facility, starting from the facility with highest ratio, i.e., the first facility in $\omega$, towards the facility with lowest ratio, i.e., the last facility in $\omega$. Every time a facility is considered, the follower allocates as many destructive resources as possible to this facility, while respecting both the individual ($y_i \leq a_i$ for all $i \in N$) and total ($\sum_{i \in N}{y_i} \leq R_f$) resource constraints. This optimal allocation for the follower is exactly how strategy $y^{\omega}$ is defined in~\eqref{eq:defyomega}, which completes the proof.
 \end{proof} \bigskip

\noindent \underline{\textbf{Proof of Lemma \ref{towardsgreedycap}}}

\begin{proof} Let $x$ be an optimal strategy and assume $x$ is neither balanced nor semi-balanced. From strategy $x$ we will construct an alternative optimal strategy $x'$ that is either balanced or semi-balanced, which proves the lemma. The outline of the proof is as follows and consists of three steps. First, we will construct an alternative strategy $x'$. Second, we will prove that strategy $x'$ is either balanced or semi-balanced. At last, we will prove that strategy $x'$ is equally good or better than strategy $x$, which implies that $x'$ is also an optimal strategy.    \bigskip

\noindent \underline{\textbf{Step 1: Constructing an alternative strategy $x'$:}} \bigskip

\noindent Let $\omega  \in \Omega^x$. For notational convenience we denote in this proof $q^{\omega}$ by $q$ and $A^{\omega}$ by $A$. In order to construct an alternative optimal strategy $x'$, we distinguish between two cases.  \bigskip

\noindent  \underline{Case 1:} Assume $R_l < \sum_{i \in N}{\frac{a_{i}p_{q}x_{q}}{a_qp_i}}$. We will construct a strategy $x'$ that is either balanced or semi-balanced at the level of $\frac{p_{q}x_{q}}{a_{q}}$. In particular, in strategy $x'$, all facilities in $A$ get resources allocated in a balanced way at level $\frac{p_{q}x_{q}}{a_{q}}$ and the facilities in $N \backslash A$ get resources  allocated in a greedily way based on their production rate, while respecting  the leader's resource constraint. Note that the facilities in $N \backslash A$ are in non-increasing order with respect to their production rates ($p_i \geq p_j$ for all $i,j \in N$ with $i<j$), which implies that we allocate these resources facility by facility (i.e., starting from the facility with lowest index towards the facility with highest index in $N \backslash A$) in a balanced way up to the level of $\frac{p_{q}x_{q}}{a_{q}}$.

Following this procedure, we define facility $k \in N \backslash A $ as the highest indexed facility that gets a positive amount of resources allocated. That means, we define facility $k \in N \backslash A $ as follows
$$k = \min \left \{j \in N \backslash A~\Bigg|~
\sum_{\substack{i \in N \backslash A: \\ i \leq j}}{\frac{a_{i}p_{q}x_{q}}{a_qp_i} } \geq R_l -  \sum_{i \in A}{\frac{a_{i}p_{q}x_{q}}{a_qp_i}} \right  \}.$$
Note that facility $k$ exists due to the assumption for Case 1. Next, define
$$B=\{i \in N \backslash A~|~i<k\}.$$
Observe that $B$ might be the empty set. Now, we are ready to introduce our strategy $x'$ formally. Strategy $x'$ is defined as
\begin{equation*}
x'_i = \begin{cases}
{\displaystyle \frac{a_{i}p_{q}x_{q}}{a_qp_i}} & \text{ if } i \in A \cup B, \\
{\displaystyle R_l -  \sum_{j \in A \cup B}{\frac{a_{j}p_{q}x_{q}}{a_qp_j}}  } & \text{ if } i=k, \\
0 & \text{ otherwise}. \end{cases}
\end{equation*}

\noindent  \underline{Case 2:}  Assume $R_l \geq \sum_{i \in N}{\frac{a_{i}p_{q}x_{q}}{a_qp_i}}$. In this case, we will construct a strategy $x'$ that is balanced at level $\frac{R_l}{\sum_{i \in N} \frac{a_i}{p_i}}$. In doing so, we define our strategy $x'$ as
\begin{equation}
\label{eq:hulp8}
x'_i = \frac{a_iR_l}{p_i \sum_{j \in N}{\frac{a_j}{p_j}}},
\end{equation}
for all $i \in N$.\bigskip

\noindent \underline{\textbf{Step 2: Proving that $x'$ is either balanced or semi-balanced:}} \bigskip

\noindent We distinguish between the previously defined two cases and we will formally prove for each case that $x'$ is either balanced or semi-balanced.   \bigskip

\noindent  \underline{Case 1:} First, observe that, based on the definition of facility $k$ and set $B$, we have
$$R_l \leq \sum_{i \in A} \frac{a_ip_qx_q}{a_qp_i} + \sum_{\substack{i \in N \backslash A: \\ i \leq k}}{\frac{a_{i}p_{q}x_{q}}{a_qp_i} } = \sum_{i \in A} \frac{a_ip_qx_q}{a_qp_i} + \sum_{i \in B} \frac{a_ip_qx_q}{a_qp_i} + \frac{a_kp_qx_q}{a_qp_k}=\sum_{i \in A \cup B \cup \{k\}} \frac{a_ip_qx_q}{a_qp_i}.$$
We will prove that $x'$ is balanced in case of an equality (i.e., $R_l=\sum_{i \in A \cup B \cup \{k\}} \frac{a_ip_qx_q}{a_qp_i}$) and semi-balanced in case of a strict inequality (i.e., $ R_l < \sum_{i \in A \cup B \cup \{k\}} \frac{a_ip_qx_q}{a_qp_i}$) in the above equation. \bigskip

\noindent  \underline{Case 1a:} Assume $R_l = \sum_{i \in A \cup B \cup \{k\}} \frac{a_ip_qx_q}{a_qp_i}$. So,  $x_k' = R_l - \sum_{j \in A \cup B} \frac{a_jp_qx_q}{a_qp_j} = \frac{a_kp_qx_q}{a_qp_k}$. We will prove that $x'$ is balanced. For this, we claim that
\begin{align*}
S^{x'}&=(A \cup B \cup \{k\}).
\end{align*}
We will now prove that $x'$ satisfies the three conditions of a balanced strategy.
\begin{enumerate}[leftmargin=0.8cm]
\item[(i)] From the definition of $x'$ it follows that
$$ \sum_{i \in N}{x'_i} = \sum_{i \in A \cup B}{x'_i} + x'_k = \sum_{i \in A \cup B}{\frac{a_{i}p_{q}x_{q}}{a_qp_i}}+ R_l -  \sum_{j \in A \cup B}{\frac{a_{j}p_{q}x_{q}}{a_qp_j}} =R_l.$$
\item[(ii)] Let $i \in S^{x'}$. We will first show  $x_q > 0$. Suppose, for the sake of contradiction, that $x_q =0$. Then, due to the assumption for Case 1, we have $R_l < \sum_{i \in N} \frac{a_ip_qx_q}{a_qp_i} = 0$, which contradicts $R_l \in \mathbb{R}_{<0}$ and thus $x_q>0$.Therefore, in combination with the definition of $x'$, we obtain $x'_i=\frac{a_ip_qx_q}{a_qp_i}>0$.
\item[(iii)] From the definition of $x'$ it follows that, for all $i,j \in S^{x'}$, we have
$\frac{p_ix'_i}{a_i}= \frac{p_qx_q}{a_q}=\frac{p_jx'_j}{a_j}.$
\item[(iv)] From the definition of $x'$ it follows that $x'_i=0$ for all $i \in N \backslash S^{x'}$.
\end{enumerate} \bigskip

\noindent  \underline{Case 1b:} Assume $R_l < \sum_{i \in A \cup B \cup \{k\}} \frac{a_ip_qx_q}{a_qp_i}$.
We will prove that $x'$ is semi-balanced.  For this, we claim that
\begin{align*}
S^{x'}&=A \cup B,  \\
r^{x'}&=k.
\end{align*}
We will now prove that $x'$ satisfies the four conditions of a semi-balanced strategy.
\begin{enumerate}[leftmargin=0.8cm]
\item[(i)] In exactly the same way as in Case 1a, we can prove that this first condition holds. \item[(ii)] In exactly the same way as in Case 1a, we can prove that this second condition holds.
\item[(iii)] In exactly the same way as in Case 1a, we can prove that this third condition holds.
\item[(iv)] From the definition of facility $k$ it follows that $\sum_{i \in A \cup B} \frac{a_ip_qx_q}{a_qp_i} < R_l$. As a consequence,
$$x'_k=R_l -  \sum_{j \in A \cup B}{\frac{a_{j}p_{q}x_{q}}{a_qp_j}} >0.$$
Next, we have
$$\sum_{i \in A \cup B \cup \{k\}} \frac{a_ip_qx_q}{a_qp_i} > R_l = \sum_{i \in N} x'_i = \sum_{i \in A \cup B} \frac{a_ip_qx_q}{a_qp_i} + x_k',$$
where the inequality follows from the assumption for Case 1b and the last equality follows from the construction of $x'$. Hence, $x'_k < \frac{a_kp_qx_q}{a_qp_k}$ and thus, for all $i \in S^{x'}$, we have
$$\frac{p_kx'_k}{a_k}< \frac{p_qx_q}{a_q}=\frac{p_ix'_i}{a_i}.$$
\item[(v)] From the definition of $x'$ it follows that $x'_i=0$ for all $i \in N \backslash (S^{x'} \cup r^{x'})$.
\end{enumerate} \bigskip

\noindent  \underline{Case 2:} We will prove that $x'$ is balanced.  For this, we claim that
\begin{align*}
S^{x'}&=N.
\end{align*}
We will now prove that $x'$ satisfies the three conditions of a balanced strategy.
\begin{enumerate}[leftmargin=0.8cm]
\item[(i)] From the definition of $x'$ it follows that
$$\sum_{i \in N}{x'_i}= \sum_{i \in N}{\frac{a_iR_l}{p_i \sum_{j \in N}{\frac{a_j}{p_j}}}}=R_l.$$
\item[(ii)] Let $i \in S^{x'}$. From the definition of $x'$ it follows that $x'_i>0$.
\item[(iii)] From the definition of $x'$ it follows that, for all $i,l \in S^{x'}$, we have
$$\frac{p_ix_i}{a_i}= \frac{R_l}{\sum_{j \in N}{\frac{a_j}{p_j}}}=\frac{p_lx_l}{a_l}.$$
\item[(iv)] Note that $N \backslash S^{x'} = \emptyset$ and thus this last condition automatically holds.
\end{enumerate}

\bigskip

\noindent \underline{\textbf{Step 3: Proving that $x'$ is also an optimal strategy:}} \bigskip

\noindent  We distinguish again between the previously defined  two cases and we will prove for each case that strategy $x'$ is equally good or better than strategy $x$. This implies that $x'$ is also an optimal strategy. For both cases we use that for every $\omega \in \Omega^x$ we have
\begin{equation}
\label{eq:aipixi>=aqomegapx>=ajpjxj}
 \frac{p_ix_i}{a_i} \geq \frac{p_{q^{\omega}}x_{q^{\omega}}}{a_{q^{\omega}}} \geq \frac{p_jx_j}{a_j},
\end{equation}
for all $i \in A^{\omega}, j \in N \backslash A^{\omega}$. These inequalities follow  from the definitions of $\Omega^{x}$, $A^{\omega}$ and $q^{\omega}$. \bigskip

\noindent \underline{Case 1:} By the construction of $x'$ we have
\begin{equation}
\label{eq:sumix'i>=sumxi}
\sum_{i \in N \backslash A}{x'_i}
 = R_l - \sum_{i \in A}{x'_i}
=R_l - \sum_{i \in A}{\frac{a_{i}p_{q}x_{q}}{a_qp_i}}
\overset{\text{\eqref{eq:aipixi>=aqomegapx>=ajpjxj}}}{\geq}   R_l - \sum_{i \in A}{x_i}
\geq \sum_{i \in N \backslash A}{x_i},
\end{equation}
where the last inequality follows from the fact $x \in \mathscr{X}$. As a consequence,
\begin{align}
\sum_{i \in N \backslash A}{x'_ip_i}=&\max_{{ (z_i)_{i \in N \backslash A}}}            & \sum_{i \in N \backslash A}{z_ip_i}&& \nonumber \\
&\text{~~~~s.t.}               & \sum_{i \in N \backslash A}{z_i}             & \leq \sum_{i \in N \backslash A}{x'_i} & \nonumber \\
&     &  z_i &\geq 0 & \quad \forall i \in N \backslash A \nonumber \\
&                & z_i              & \leq \frac{a_{i}p_{q}x_{q}}{a_qp_i} &\quad \forall i \in N \backslash A \nonumber   \\
\overset{\text{\eqref{eq:aipixi>=aqomegapx>=ajpjxj},\eqref{eq:sumix'i>=sumxi}}}{\geq}&  \sum_{i \in N \backslash A}{x_ip_i}, \label{eq:sumix'ipi>=sumxipi}
\end{align}
where the equality follows from the fact that the construction of $x'$ is such that the facilities in $N \backslash A$ get resources  allocated in a greedily way based on their production rate, while respecting  the leader's resource constraint and the facility's capacity constraint. The inequality follows from the fact that $(x_i)_{i \in N}$ lies in the feasible facility of the corresponding LP-problem, due to~\eqref{eq:sumix'i>=sumxi} (implying the first constraint), due to~\eqref{eq:aipixi>=aqomegapx>=ajpjxj} (implying the third constraint), and due to the fact $x \in  \mathscr{X}$ and thus $x_i \geq 0$ for all $i \in N \backslash A$ (implying the second constraint).

Now, observe that, by construction of $x'$, we have, for $j \in B$,
\begin{equation} \label{eq:hulp1} \frac{p_jx'_j}{a_j} = \frac{p_{q}x_{q}}{a_{q}}.
\end{equation}
Next, as shown in Case 1a and Case 1b of Step 2, we have
\begin{equation} \label{eq:hulp2} x_k' \leq \frac{a_kp_qx_q}{a_qp_k}. \end{equation}
At last, by construction of $x'$, we have, for $j \in N \backslash (A \cup B \cup \{k\})$,
\begin{equation} \label{eq:hulp3}  x_j' = 0 \leq \frac{a_jp_qx_q}{a_qp_j}. \end{equation}
Hence, from \eqref{eq:hulp1}, \eqref{eq:hulp2} and \eqref{eq:hulp3} it follows that, for all $j \in N \backslash A$, we have
\begin{equation}
\label{eq:hulp4}
    \frac{p_qx_q}{a_q} \geq \frac{p_jx'_j}{a_j}.
\end{equation}
Note that, by the construction of $x'$, we also have
\begin{equation}
\label{eq:x'j*=xj*}
x'_{q} = \frac{a_qp_qx_q}{a_qp_q} = x_{q}.
\end{equation}
Hence, for all $i \in A, j \in N \backslash A$, we have
\begin{equation}
\label{eq:aipix'i>=aqomegapx'>=ajpjx'j}
\frac{p_ix'_i}{a_i} = \frac{p_{q}x_{q}}{a_{q}}  \overset{\text{\eqref{eq:x'j*=xj*}}}{=}   \frac{p_{q}x'_{q}}{a_{q}} \overset{\text{\eqref{eq:x'j*=xj*}}}{=}   \frac{p_{q}x_{q}}{a_{q}}  \overset{\text{\eqref{eq:hulp4}}}{\geq}  \frac{p_jx'_j}{a_j},
\end{equation}
where the first equality follows from the construction of $x'$. It follows from~\eqref{eq:aipix'i>=aqomegapx'>=ajpjx'j} that there exists a permutation $\omega' \in \Omega^{x'}$ such that $q^{\omega'}=q$ and $A^{\omega'}=A$. As a consequence, also $y^{\omega'}=y^{\omega}$. Hence,
{\allowdisplaybreaks
\begin{align*}
 \mathscr{P}(x') & \overset{\text{\eqref{eq:defPxyomega}}}{=}  x'_{q^{\omega'}}p_{q^{\omega'}} \left (1-\frac{y^{\omega'}_{q^{\omega'}}}{a_{q^{\omega'}}} \right ) + \sum_{i \in N \backslash A^{\omega'}}{x'_ip_i} \\
& = x'_{q}p_{q} \left (1-\frac{y^{\omega}_{q}}{a_{q}} \right ) + \sum_{i \in N \backslash A}{x'_ip_i} \\
& \stackrel{\mathclap{\text{\eqref{eq:x'j*=xj*},\eqref{eq:sumix'ipi>=sumxipi}}}}{\geq}   x_{q}p_{q} \left (1-\frac{y^{\omega}_{q}}{a_{q}} \right ) + \sum_{i \in N \backslash A}{x_ip_i} \\
 &  \overset{\text{\eqref{eq:defPxyomega}}}{=} \mathscr{P}(x).
\end{align*}}Hence, strategy $x'$ is equally good or better than strategy $x$ and thus $x'$ is also an optimal strategy. \bigskip

\noindent \underline{Case 2:} We will first define an alternative strategy, called strategy $\hat{x}$. In this strategy, every facility in $N$ gets resources allocated in a balanced way, at level $\frac{p_qx_q}{a_q}$. We will prove that strategy $\hat{x}$ is equally good or better than strategy $x$. Thereafter, we will prove that strategy $x'$ is equally good or better than strategy $\hat{x}$, implying that strategy $x'$ is also equally good or better than strategy $x$. As a result, $x'$ is also an optimal strategy.

Define strategy $\hat{x}$ as follows
\begin{equation}
    \label{eq:hulp7}
\hat{x}_i=\frac{a_{i}p_{q}x_{q}}{a_qp_i},
\end{equation}
for all $i \in N$. Now, observe that, for all $i \in (N \backslash A) \cup \{q\}$, it follows that
\begin{equation}
\label{eq:hulp6}
x_i \overset{\text{\eqref{eq:aipixi>=aqomegapx>=ajpjxj}}}{\leq}  \frac{a_{i}p_{q}x_{q}}{a_qp_i}
\overset{\text{\eqref{eq:hulp7}}}{=}  \hat{x}_i,
\end{equation}
where the inequality is an equality for $i=q$. Moreover, observe that, for all $i,j \in N$, we have
$$\frac{p_i\hat{x}_i}{a_i} \overset{\text{\eqref{eq:hulp7}}}{=} \frac{p_q x_q}{a_q} \overset{\text{\eqref{eq:hulp7}}}{=} \frac{p_j\hat{x}_j}{a_j},$$
and thus $\omega \in \Omega^{\hat{x}}$. From this we conclude
{\allowdisplaybreaks
\begin{equation}
\label{eq:Phatx>=P(x)}
 \mathscr{P}(\hat{x})  \overset{\text{\eqref{eq:defPxyomega}}}{=}  \hat{x}_{q}p_{q} \left (1-\frac{y^{\omega}_{q}}{a_{q}} \right ) + \sum_{i \in N \backslash A}{\hat{x}_ip_i}  \overset{\text{\eqref{eq:hulp6}}}{\geq}    x_{q}p_{q} \left (1-\frac{y^{\omega}_{q}}{a_{q}} \right ) + \sum_{i \in N \backslash A}{x_ip_i}  \overset{\text{\eqref{eq:defPxyomega}}}{=} \mathscr{P}(x).
\end{equation}}

\noindent Now, we will focus on strategy $x'$. Recall from~\eqref{eq:hulp8} that in strategy $x'$ every facility in $N$ gets resources allocated in balanced way, at level $\frac{R_l}{\sum_{i \in N} \frac{a_i}{p_i}}$. First, observe that, for all $i \in N$, we have
\begin{equation}
\label{eq:hulpje}
x'_i \overset{\text{\eqref{eq:hulp8}}}{=} \frac{a_iR_l}{p_i \sum_{j \in N}{\frac{a_j}{p_j}}}
\geq \frac{a_i \sum_{j \in N}{\frac{a_{j}p_{q}x_{q}}{a_qp_j}}}{p_i \sum_{j \in N}{\frac{a_j}{p_j}}}
= \frac{a_{i}p_{q}x_{q}}{a_qp_i} \overset{\text{\eqref{eq:hulp7}}}{=} \hat{x}_i,\end{equation}
where the inequality follows from the assumption for Case 2.
Moreover, observe that, for all $i,j \in N$, we have $$\frac{p_ix'_i}{a_i} \overset{\text{\eqref{eq:hulp8}}}{=} \frac{R_l}{\sum_{k \in S} \frac{a_k}{p_k}} \overset{\text{\eqref{eq:hulp8}}}{=} \frac{p_jx'_j}{a_j},$$
and thus $\omega \in \Omega^{x'}$. From this we conclude
{\allowdisplaybreaks
\begin{equation*}
 \mathscr{P}(x')  \overset{\text{\eqref{eq:defPxyomega}}}{=}  x'_{q}p_{q} \left (1-\frac{y^{\omega}_{q}}{a_{q}} \right ) + \sum_{i \in N \backslash A}{x'_ip_i}  \overset{\text{\eqref{eq:hulpje}}}{\geq}    \hat{x}_{q}p_{q} \left (1-\frac{y^{\omega}_{q}}{a_{q}} \right ) + \sum_{i \in N \backslash A}{\hat{x}_ip_i}  \overset{\text{\eqref{eq:defPxyomega}}}{=} \mathscr{P}(\hat{x}),
\end{equation*}}
where the last equality also follows from $\omega  \in \Omega ^{\hat{x}}$. This implies together with~\eqref{eq:Phatx>=P(x)} that
$$
 \mathscr{P}(x')  \geq  \mathscr{P}(x),
 $$
and thus $x'$ is also an optimal strategy.
\end{proof} \bigskip

\noindent \underline{\textbf{Proof of Lemma \ref{semi-balanced}}}

\begin{proof}
Let $\omega \in \Omega^x$. For notational convenience we denote in this proof $S^x$ by $S$, $r^{x}$ by $r$, $q^{\omega}$ by $q$ and $A^{\omega}$ by $A$. From the definition of a semi-balanced strategy it follows immediately that, for all $i \in N \backslash (S \cup \{r\})$, we have
\begin{equation}
\label{eq:xi=0foralliinNSxcupi*}
x_i = 0.
\end{equation}
Next, let $i \in S$. It follows from the definition of a semi-balanced strategy that, for all $j \in S$, we have $\frac{p_ix_i}{a_i}=\frac{p_jx_j}{a_j}$ and thus $x_j=\frac{a_jp_ix_i}{a_ip_j}$.  As a consequence,
$$R_l=\sum_{j \in N}{x_j} \overset{\text{\eqref{eq:xi=0foralliinNSxcupi*}}}{=} x_{r}+\sum_{j \in S}{x_j}=x_{r}+\sum_{j \in S}{\frac{a_jp_ix_i}{a_ip_j}} = x_{r}+\frac{p_ix_i}{a_i}\sum_{j \in S}{\frac{a_j}{p_j}},$$
where the first equality follows from the definition of a semi-balanced strategy. Hence, for all $i \in S$, we have
\begin{equation}
\label{eq:xiforalliinSx}
x_i = \frac{a_i (R_l - x_{r})}{p_i \sum_{j \in S} \frac{a_j}{p_j}},
\end{equation}
which concludes the first part of the lemma.

 For the second part of the lemma, we distinguish from now on between five different cases. \bigskip

\noindent \underline{Case 1:} Assume $R_f < \sum_{i \in S} a_i$. Then, $A \subseteq S$ and thus
\begin{align*} \mathscr{P}(x)
 & \overset{\text{\eqref{eq:defPxyomega}}}{=}  x_{q}p_{q}\left (1-\frac{y^{\omega}_{q}}{a_{q}} \right ) + \sum_{i \in N \backslash A}{x_ip_i} \\
 &=  -\frac{x_{q}p_{q}}{a_{q}}y^{\omega}_{q}  + \sum_{i \in (N \backslash A) \cup \{q\}}{x_ip_i} \\
 &=  -\frac{x_{q}p_{q}}{a_{q}}y^{\omega}_{q}  + \sum_{i \in (S \backslash A) \cup \{q\}}{x_ip_i} + \sum_{i \in N \backslash S}{x_ip_i} \\
  &=  -\frac{x_{q}p_{q}}{a_{q}}y^{\omega}_{q}  + \sum_{i \in (S \backslash A) \cup \{q\}}{x_ip_i} +x_{r}p_{r}  + \sum_{i \in N \backslash (S \cup \{r\})}{x_ip_i} \\
 &\stackrel{\mathclap{\text{\eqref{eq:defyomega},\eqref{eq:xi=0foralliinNSxcupi*},\eqref{eq:xiforalliinSx}}}}{=} \hspace*{1.2em}   -\frac{(R_l-x_{r})}{\sum_{j \in S} \frac{a_j}{p_j}} \left (R_f - \sum_{i \in A \backslash \{q\}}{a_i}  \right )  + \sum_{i \in (S \backslash A) \cup \{q\}}{\frac{a_i(R_l-x_{r})}{ \sum_{j \in S} \frac{a_j}{p_j}}} +x_{r}p_{r} \\
&=  \frac{\left(\sum_{i \in S} a_i - R_f\right) \left(R_l-x_{r}\right)}{ \sum_{i \in S} \frac{a_i}{p_i}} + p_{r} x_{r}, \\
&=\overline{p}(S)(R_l-x_{r}) + p_{r} x_{r},
\end{align*}
where the second equality holds because $q \in A$.  The third equality follows from the fact that $N \backslash A$ can be partitioned in the sets $S \backslash A$ and $N \backslash S$, due to the fact that $A \subseteq S$. The fourth equality holds because $r \not \in S$. For the penultimate equality we use $q \in A \subseteq S$ and thus $((S \backslash A) \cup \{q\}) \cup (A \backslash \{q\}) = S$. For the last equality we use $R_l < \sum_{i \in S}a_i $ (the assumption for Case 1) and thus $\overline{p}(S) = \frac{\left(\sum_{i \in S} a_i - R_f\right) \left(R_l-x_{r}\right)}{ \sum_{i \in S} \frac{a_i}{p_i}}$.\bigskip

\noindent \underline{Case 2:} Assume $R_f = \sum_{i \in S} a_i$ and thus $\overline{p}(S)=0$. Then, we have $A=S$ and $y^{\omega}_{q}=a_{q}$. As a consequence,
\begin{align*} \mathscr{P}(x)
  &\overset{\text{\eqref{eq:defPxyomega}}}{=}  x_{q}p_{q}\left (1-\frac{y^{\omega}_{q}}{a_{q}} \right ) + \sum_{i \in N \backslash A}{x_ip_i} \\
  &  =  x_{q}p_{q}\left (1-\frac{a_{q}}{a_{q}} \right ) + \sum_{i \in N \backslash S}{x_i p_i} \\
  &  \overset{\text{\eqref{eq:xi=0foralliinNSxcupi*}}}{=} x_{r}p_{r}.
\end{align*}

\noindent \underline{Case 3:} Assume $\sum_{i \in S} a_i < R_f < \sum_{i \in S \cup \{r\}} a_i$. Then, we have $A=S \cup \{r\}$ and $q=r$. As a consequence,
\begin{align*} \mathscr{P}(x)
  &\overset{\text{\eqref{eq:defPxyomega}}}{=}  x_{q}p_{q}\left (1-\frac{y^{\omega}_{q}}{a_{q}} \right ) + \sum_{i \in N \backslash A}{x_ip_i} \\
  &=  x_{r}p_{r}\left (1-\frac{y^{\omega}_{q}}{a_{r}} \right ) + \sum_{i \in N \backslash (S \cup \{r\})}{x_ip_i} \\
  & \stackrel{\mathclap{\text{\eqref{eq:defyomega},\eqref{eq:xi=0foralliinNSxcupi*}}}}{=} \hspace*{0.3em}  x_{r}p_{r} \left (1-\frac{R_f - \sum_{i \in S}{a_i}}{a_{r}} \right )  \\
  &= \frac{p_{r}x_{r}}{a_{r}} \left( \sum_{i \in S \cup \{r\}} a_i-R_f\right).
\end{align*}

\noindent \underline{Case 4:} Assume $R_f = \sum_{i \in S \cup \{r\}} a_i$. Then, $A=S \cup \{r\}$ and $y^{\omega}_{q}=a_{q}$. As a consequence,
$$ \mathscr{P}(x)
  \overset{\text{\eqref{eq:defPxyomega}}}{=}  x_{q}p_{q}\left (1-\frac{y^{\omega}_{q}}{a_{q}} \right ) + \sum_{i \in N \backslash A}{x_ip_i}  =  x_{q}p_{q}\left (1-\frac{a_{q}}{a_{q}} \right ) + \sum_{i \in N \backslash (S \cup \{r\})}{x_i  p_i}
 \overset{\text{\eqref{eq:xi=0foralliinNSxcupi*}}}{=}0.
$$

\noindent \underline{Case 5:} Assume $R_f > \sum_{i \in S \cup \{r\}} a_i$. Then, $S \cup \{r\} \subseteq A \backslash \{q\}$ and thus $(N \backslash A) \cup \{q\} = N \backslash (A \backslash \{q\})$ $\subseteq N \backslash (S \cup \{r\})$. Hence, it follows from~\eqref{eq:xi=0foralliinNSxcupi*} that $x_i=0$ for all $i \in (N \backslash A) \cup \{q\}$. As a consequence,
\begin{equation*} \mathscr{P}(x)
  \overset{\text{\eqref{eq:defPxyomega}}}{=}  x_{q}p_{q}\left (1-\frac{y^{\omega}_{q}}{a_{q}} \right ) + \sum_{i \in N \backslash A}{x_ip_i}   =0. \qedhere
\end{equation*}

\end{proof} \bigskip

\noindent \underline{\textbf{Proof of Lemma \ref{allocationbalancedx+balanced}}}

\begin{proof}
This lemma follows directly from the proof of Lemma~\ref{semi-balanced} by applying the definition of composed net production rate (Definition~\ref{def:avgrate}) and taking $x_{r}=0$. (Note that this is possible because that proof doesn't require $x_{r} > 0$ or $x_{r} \neq 0$ in any of the arguments.)
\end{proof} \bigskip

\noindent \underline{\textbf{Proof of Lemma \ref{t:nofractions}}}

\begin{proof} Let $x$ be an optimal strategy that is either balanced or semi-balanced. Note that such a strategy always exists due to Lemma~\ref{towardsgreedycap}. If $x$ is balanced, then the theorem holds immediately. Therefore, we assume $x$ is semi-balanced.  From strategy $x$ we will construct two alternative balanced strategies, called strategy $x'$ and $x''$. Next, we will prove that at least one of them is equally good or better than strategy $x$, which implies that this strategy is also an optimal strategy.

For notational convenience we denote in this proof $S^x$ by $S$ and $r^{x}$ by $r$. Next, let $x'$ and $x''$ be the two unique balanced strategies  such that $S^{x'} = S$ and  $S^{x''} = S \cup \{r\}$. We distinguish between four cases and for each case we show that either $\mathscr{P}(x') \geq \mathscr{P}(x)$ or $\mathscr{P}(x'') \geq \mathscr{P}(x)$ holds. \bigskip

\noindent \underline{Case 1:} Assume $R_f \geq \sum_{i \in S \cup \{r\}}{a_i}$. Then, by Corollary~\ref{cor:allocationbalancedx+balancedNIEUW} and Lemma~\ref{semi-balanced} we have
$$\mathscr{P}(x'') =\overline{p}(S \cup \{r\})R_l = 0 \cdot R_l = 0 = \mathscr{P}(x).$$

\noindent \underline{Case 2:} Assume  $\sum_{i \in S}{a_i} < R_f < \sum_{i \in S \cup \{r\}}{a_i}$. From the definition of a semi-balanced strategy it follows immediately that $\frac{p_ix_i}{a_i} > \frac{p_{r}x_{r}}{a_{r}}$  for all $i \in S$. Hence,
$$R_l = \sum_{i \in S \cup \{r\}}{x_i} > \sum_{i \in S \cup \{r\}}{\frac{a_ip_{r}x_{r}}{a_{r}p_i}} = \frac{p_{r}x_{r}}{a_{r}}\sum_{i \in S \cup \{r\}}{\frac{a_i}{p_i}}.$$
As a consequence, by Corollary~\ref{cor:allocationbalancedx+balancedNIEUW} and Lemma~\ref{semi-balanced}, we have
$$\mathscr{P}(x'')
=\overline{p}(S \cup \{r\})R_l
= \frac{R_l}{\sum_{i \in S \cup \{r\}}{\frac{a_i}{p_i}}}\left( \sum_{i \in S \cup \{r\}} a_i-R_f\right)
> \frac{p_{r}x_{r}}{a_{r}} \left( \sum_{i \in S \cup \{r\}} a_i-R_f\right)
= \mathscr{P}(x).
$$

\noindent \underline{Case 3:} Assume $R_f \leq  \sum_{i \in S}{a_i}$ and $ p_{r} \leq  \overline{p}(S) $. Then, by Corollary~\ref{cor:allocationbalancedx+balancedNIEUW} and Lemma~\ref{semi-balanced}, it follows  immediately  that
$$\mathscr{P}(x')
= \overline{p}(S)R_l
\geq \overline{p}(S) (R_l-x_{r}) + p_{r} x_{r}
= \mathscr{P}(x).
$$

\noindent \underline{Case 4:} Assume $R_f \leq  \sum_{i \in S}{a_i}$ and $ p_{r} >  \overline{p}(S)$. From the definition of a semi-balanced strategy together with Lemma~\ref{semi-balanced}, it follows that, for all $i \in S$, we have
\begin{equation}
\label{eq:tweedehulpvoorCase4}
\frac{p_{r}x_{r}}{a_{r}}
<\frac{p_ix_i}{a_i}
= \frac{R_l - x_{r}}{ \sum_{j \in S} \frac{a_j}{p_j}}.
\end{equation}
Next, note that, by Corollary~\ref{cor:allocationbalancedx+balancedNIEUW} and Lemma~\ref{semi-balanced}, we have
\begin{align*}
&\mathscr{P}(x'') - \mathscr{P}(x)
=\overline{p} \left (S \cup \{r\} \right )R_l - \overline{p} (S)(R_l-x_r) + p_{r} x_{r} \\
&= \frac{\left(\sum_{i \in S \cup \{r\}} a_i - R_f\right) R_l}{ \sum_{i \in S \cup \{r\}} \frac{a_i}{p_i}}
- \left ( \frac{\left(\sum_{i \in S} a_i - R_f\right) \left(R_l-x_{r}\right)}{ \sum_{i \in S} \frac{a_i}{p_i}} + p_{r} x_{r} \right ) \\
&= \frac{{\displaystyle \sum_{i \in S} \frac{a_i}{p_i} \left(\sum_{i \in S \cup \{r\}} a_i - R_f\right) R_l
 -  \sum_{i \in S \cup \{r\}} \frac{a_i}{p_i} \left(\sum_{i \in S} a_i - R_f\right) \left(R_l-x_{r}\right)
 - p_{r} x_{r} \sum_{i \in S} \frac{a_i}{p_i} \sum_{i \in S \cup \{r\}} \frac{a_i}{p_i}}}
 {\sum_{i \in S} \frac{a_i}{p_i} \sum_{i \in S \cup \{r\}} \frac{a_i}{p_i}}.
\end{align*}
We will prove that the numerator of the latter expression is positive.
\begin{align*}
& \sum_{i \in S} \frac{a_i}{p_i} \left(\sum_{i \in S \cup \{r\}} a_i - R_f\right) R_l
 -  \sum_{i \in S \cup \{r\}} \frac{a_i}{p_i} \left(\sum_{i \in S} a_i - R_f\right) \left(R_l-x_{r}\right)
 - p_{r} x_{r} \sum_{i \in S} \frac{a_i}{p_i} \sum_{i \in S \cup \{r\}} \frac{a_i}{p_i} \\
 &= \sum_{i \in S} \frac{a_i}{p_i} \left(\sum_{i \in S} a_i - R_f\right) R_l + a_rR_l\sum_{i \in S} \frac{a_i}{p_i}
 -  \sum_{i \in S } \frac{a_i}{p_i} \left(\sum_{i \in S} a_i - R_f\right) \left(R_l-x_{r}\right) \\
 &- \frac{a_r}{p_r}\left(\sum_{i \in S} a_i - R_f\right) \left(R_l-x_{r}\right)
 - p_{r} x_{r} \sum_{i \in S} \frac{a_i}{p_i} \sum_{i \in S } \frac{a_i}{p_i}
  - a_{r} x_{r} \sum_{i \in S} \frac{a_i}{p_i} \\
 &= a_{r}R_l\sum_{i \in S } \frac{a_i}{p_i}
  + x_{r}\sum_{i \in S} \frac{a_i}{p_i} \left(\sum_{i \in S } a_i - R_f\right)
  -  \frac{a_{r}}{p_{r}} \left(\sum_{i \in S} a_i - R_f\right) \left(R_l-x_{r}\right)\\
  &- p_{r} x_{r} \sum_{i \in S} \frac{a_i}{p_i} \sum_{i \in S } \frac{a_i}{p_i}
  - a_{r} x_{r} \sum_{i \in S} \frac{a_i}{p_i}\\
 &=  \frac{a_{r}}{p_{r}} \left(R_l - x_{r}  \right )p_{r} \sum_{i \in S}{\frac{a_i}{p_i}}
 + \frac{a_{r}}{p_{r}} \left(R_l - x_{r}  \right ) \left ( R_f - \sum_{i \in S}{a_i}  \right )
 -  x_{r} \sum_{i \in S} \frac{a_i}{p_i}p_{r} \sum_{i \in S}{\frac{a_i}{p_i}} \\
 &-  x_{r} \sum_{i \in S} \frac{a_i}{p_i}  \left ( R_f - \sum_{i \in S}{a_i}  \right )
  \\
   &=\left( \frac{a_{r}}{p_{r}} \left(R_l - x_{r}  \right ) -  x_{r} \sum_{i \in S} \frac{a_i}{p_i} \right )
\left ( p_{r} \sum_{i \in S}{\frac{a_i}{p_i}} + R_f - \sum_{i \in S}{a_i}  \right ) \\
&>0,
\end{align*}
where the inequality follows from~\eqref{eq:tweedehulpvoorCase4} together with the assumption $ p_{r} >  \overline{p}(S)$ for Case 4. Hence, we can conclude
\begin{equation*}
\mathscr{P}(x'') > \mathscr{P}(x).\qedhere
\end{equation*}
\end{proof} \bigskip
 \newpage
\noindent \underline{\textbf{Proof of Lemma \ref{veiligtoevoegen}}}

\begin{proof}
For notational convenience we denote in this proof $S^x$ by $S$ and $S^{x'}$ by $T$. We distinguish between three cases. \bigskip

\noindent \underline{Case 1:} Assume $\sum_{j \in S }{a_j} < \sum_{j \in S \cup T}{a_j} \leq R_f$. Then, $\overline{p}(S)=0$ and $\overline{p}(S \cup T)=0$, so all four statements  follow immediately. \bigskip

\noindent \underline{Case 2:} Assume $\sum_{j \in S }{a_j} \leq R_f < \sum_{j \in S \cup T}{a_j} $. Then, $\overline{p}(S)=0$ and $\overline{p}(S \cup T) = \frac{\sum_{i \in S \cup T}a_i - R_f}{\sum_{i \in S \cup T} \frac{a_i}{p_i}}>0$, so statement \emph{(i)} and \emph{(iii)} follow immediately. Moreover, the if condition of statement \emph{(ii)} in Lemma~\ref{veiligtoevoegen},  $p_i \leq \overline{p}(S)=0$  for all $i \in T$,   is not possible, because the production rates are assumed to be strictly positive. We will now show that also the if condition of statement \emph{(iv)} in Lemma~\ref{veiligtoevoegen} is not possible. For this, suppose we have, for all $i \in T$,
$$p_i \leq \overline{p}(S \cup T) = \frac{\sum_{j \in S \cup T} a_j - R_f }{ \sum_{j \in S \cup T} \frac{a_j}{p_j}}.$$
Next, let $k \in \argmax_{i \in T}{p_i}$, then it follows that
\begin{equation}
\label{eq:lemmaveiligtoevoegenCase2(iv)}
\sum_{j \in S \cup T}{\frac{a_j}{p_j}} - \frac{1}{p_k} \left ( \sum_{j \in S \cup T}{a_j}-R_f    \right)  \leq 0.
\end{equation}
As a consequence,
\begin{align*}
\sum_{j \in S}{\frac{a_j}{p_j}} - \frac{1}{p_k} \left ( \sum_{j \in S}{a_j}-R_f    \right)
&= \sum_{j \in S \cup T}{\frac{a_j}{p_j}} - \frac{1}{p_k} \left ( \sum_{j \in S \cup T}{a_j}-R_f    \right) - \left (  \sum_{i \in T}{\frac{a_i}{p_i}} - \frac{1}{p_k}  \sum_{i \in T}{a_i}   \right ) \\
&= \sum_{j \in S \cup T}{\frac{a_j}{p_j}} - \frac{1}{p_k} \left ( \sum_{j \in S \cup T}{a_j}-R_f    \right) - \sum_{i \in T}{a_i} \left (\frac{1}{p_i}-\frac{1}{p_k}   \right )  \\
& \overset{\text{\eqref{eq:lemmaveiligtoevoegenCase2(iv)}}}{\leq} 0,
\end{align*}
where the inequality also uses $k \in \argmax_{i \in T}{p_i}$ and thus $\frac{1}{p_i}- \frac{1}{p_k} \geq 0$ for all $i \in T$. Hence, we have
$$p_k \leq  \frac{\sum_{j \in S } a_j - R_f}{ \sum_{j \in S } \frac{a_j}{p_j}} \leq 0,$$
where the last inequality follows from the assumption for Case 2. Hence, we have $p_i \leq 0$ for all $i \in T$, which is not possible, because the production rates are assumed to be strictly positive.

\bigskip

\noindent \underline{Case 3:} Assume $R_f < \sum_{j \in S }{a_j} < \sum_{j \in S \cup T}{a_j} $. We will now prove every statement separately. \bigskip

\noindent \emph{(i)} Assume $p_i \geq \overline{p}(S)$ for all $i \in T$. Then, for all $ i \in T$, we have
\begin{equation}
\label{eq:lemmaveiligtoevoegen(i)}
\sum_{j \in S}{\frac{a_j}{p_j}} - \frac{1}{p_i} \left ( \sum_{j \in S}{a_j}-R_f    \right)  \geq 0,
\end{equation}
As a consequence,
\begin{align*}
\overline{p}(S \cup T) - \overline{p}(S)
&= \frac{{\displaystyle \sum_{j \in S \cup T} a_j - R_f }}{ \sum_{j \in S \cup T} \frac{a_j}{p_j}} -  \frac{{\displaystyle \sum_{j \in S} a_j - R_f }}{ \sum_{j \in S} \frac{a_j}{p_j}} \\
&=   \frac{{\displaystyle   \sum_{j \in S} \frac{a_j}{p_j}
 \left(\sum_{j \in S \cup T} a_j - R_f\right) -
 \sum_{j \in S \cup T} \frac{a_j}{p_j}
 \left(\sum_{j \in S} a_j - R_f\right)  }}
 {\sum_{j \in S \cup T} \frac{a_j}{p_j} \sum_{j \in S} \frac{a_j}{p_j}} \\
& =  \frac{{\displaystyle \sum_{j \in S} \frac{a_j}{p_j}
 \sum_{j \in T} a_j  -
 \sum_{j \in  T} \frac{a_j}{p_j}
 \left(\sum_{j \in S} a_j - R_f\right)   }}
 {\sum_{j \in S \cup T} \frac{a_j}{p_j} \sum_{j \in S} \frac{a_j}{p_j}} \\
 & =  \frac{{\displaystyle \sum_{i \in T}{a_i}
 \left ( \sum_{j \in S}{\frac{a_j}{p_j}} - \frac{1}{p_i} \left ( \sum_{j \in S}{a_j}-R_f    \right) 	\right )}}
 {\sum_{j \in S \cup T} \frac{a_j}{p_j} \sum_{j \in S} \frac{a_j}{p_j}}  \\
 & \overset{\text{\eqref{eq:lemmaveiligtoevoegen(i)}}}{\geq} 0.
\end{align*}

\noindent \emph{(ii)} Assume $p_i \leq  \overline{p}(S)$ for all $i \in T$. In a similar way as in the proof of statement \emph{(i)} in Case 3 we can prove that $\overline{p}(S \cup T) \leq \overline{p}(S)$. \bigskip

\noindent \emph{(iii)} Assume $p_i \geq \overline{p}(S \cup T)$ for all $i \in T$. Let $k \in \argmin_{i \in T}{p_i}$, then we can prove in a similar way as in Case 2 that
$$p_k \geq  \frac{\sum_{j \in S } a_j - R_f }{ \sum_{j \in S } \frac{a_j}{p_j}} = \overline{p}(S).$$
 Hence, $p_i \geq \overline{p}(S)$ for all $i \in T$ and thus it  follows immediately from Lemma~\ref{veiligtoevoegen}\emph{(i)} that  $\overline{p}(S \cup T) \geq \overline{p}(S)$. \bigskip

\noindent \emph{(iv)} Assume $p_i \leq \overline{p}(S \cup T)$ for all $i \in T$. In a similar way as in the proof of statement \emph{(iii)} in Case 3, and using $k \in \argmax_{i \in T}{p_i}$,  we can prove that $\overline{p}(S \cup T) \leq \overline{p}(S)$.
\end{proof} \bigskip

\noindent \underline{\textbf{Proof of Lemma \ref{bijnaklaarstelling}}}

\begin{proof}
Let $x$ be an optimal strategy that is balanced. Note that such a strategy always exists due to Lemma~\ref{t:nofractions}. Assume $x$ is not seried. From strategy $x$ we will construct an alternative optimal balanced strategy $x''$ that is seried, which proves the lemma. The outline of the proof is as follows and consists of the following steps. First, we add facilities to $S^x$ for which the production rate is at least the composed net production rate $\overline{p}(S^x)$. Second, we exclude, from this new set, the facilities  for which the production rate is less than the production rate production rate $\overline{p}(S^x)$. Next, we show that the resulting new balanced strategy is equally good or better than strategy $x$, which implies that this new balanced strategy is also optimal. At last, we will prove that this new optimal balanced strategy is seried.

For notational convenience we denote in this proof $S^x$ by $S$. Next, denote the facility in $S$ with the highest index by $s$, i.e.,
$$s=\max\{i~|~i \in S\}.$$  Note that since $x$ is not seried, we have $S \subset \{1,\ldots,s\}$. Let $x'$ be the balanced strategy with
$$S^{x'}=S \cup \left \{i \in \{1,\ldots,s\} \backslash S~|~p_i \geq \overline{p}(S) \right  \}.$$
Then, by Lemma~\ref{veiligtoevoegen}\emph{(i)}, we have
\begin{equation}
\label{eq:bijnaklaarstellingx'geqx}
\overline{p}(S^{x'}) \geq \overline{p}(S).
\end{equation}

\noindent Note that if $p_i \geq \overline{p}(S)$ for all $i \in \{1,\ldots,s\} \backslash S$, then $S^{x'} = \{1,\ldots,s\}$ and thus $x'$ is a seried strategy.

On the other hand, if this is not the case, so there exists a facility $i \in \{1,\ldots,s\} \backslash S$ with $p_i < \overline{p}(S)$, then $x'$ is not seried. Therefore, we will continue to find a seried strategy as follows. Let $x''$ be the balanced strategy such that
$$S^{x''}=S^{x'} \big \backslash \left \{i \in S^{x'} ~|~p_i < \overline{p}(S) \right  \}.$$
Note that, by~\eqref{eq:bijnaklaarstellingx'geqx} we have that $p_i < \overline{p}(S)$ implies also $p_i < \overline{p}(S^{x'})$. As a consequence, by Lemma~\ref{veiligtoevoegen}\emph{(iv)}, we have $$\overline{p}(S^{x''}) \geq \overline{p}(S^{x'}) \overset{\text{\eqref{eq:bijnaklaarstellingx'geqx}}}{\geq} \overline{p}(S),$$ and thus
\begin{equation}
\label{eq:bijnaklaarstellingPxy>=Px''y''}
\mathscr{P}(x'') \geq \mathscr{P}(x).
\end{equation}
Hence, $x''$ is also an optimal balanced strategy.

We will now show that $x''$ is a seried strategy. We do this by first showing that $S^{x''}$ is non-empty. After that, we show that every facility with an index lower than index $t$, where facility $t$ is the facility with the highest index in $S^{x''}$, is also an element of $S^{x''}$.

We will prove that $S^{x''}$ is non-empty by showing that facility $1$ is an element of $S^{x''}$. Since $p_1 \geq p_i$ for all $i \in N$, we have
$$\sum_{i \in S}{\frac{a_i}{p_i}} - \frac{1}{p_1} \left ( \sum_{i \in S}{a_i}-R_f    \right) = \sum_{i \in S}{a_i\left( \frac{1}{p_i}-\frac{1}{p_1} \right )} +\frac{R_f}{p_1}  \geq 0,$$
and thus
$$p_1 \geq \frac{\sum_{i \in S} a_i - R_f}{ \sum_{i \in S} \frac{a_i}{p_i}}.$$
From this, together with the fact $p_1>0$, we can conclude $p_1 \geq \max\left\{\frac{\sum_{i \in S} a_i - R_f}{\sum_{i \in S} \frac{a_i}{p_i}},0\right\} = \overline{p}(S)$ and thus facility 1 is an element of $S^{x''}$, which implies that $S^{x''} \neq \emptyset$.

Next, let $t = \max\{i~|~i \in S^{x''}\}$. Note that $t$ exists because $S^{x''} \neq \emptyset$. Moreover, note that  $t \leq s$ by construction of $S^{x'}$ and $S^{x''}$. In order to prove that $x''$ is a seried strategy, we need to show that, for all $i \in N$ with $i<t$, we also have $i \in S^{x''}$. For this, let $i \in N$ with $i<t$. Note that, since $t \in S^{x''}$, we have $p_t \geq  \overline{p}(S)$. As a consequence, because $i<t$, we have $p_i \geq p_t$ and thus also
\begin{equation}
\label{eq:pi>=Pxyomega/Rd}
p_i \geq  \overline{p}(S).
\end{equation}
In order to prove that $i \in S^{x''}$, we first show that $i \in S^{x'}$. For this, we distinguish between two cases: $i \in S$ and $i \not \in S$. Note that if $i \in S$, then also $i \in S^{x'}$ because $S \subseteq S^{x'}$. On the other hand, if $i \not \in S$, then from~\eqref{eq:pi>=Pxyomega/Rd} together with the fact $i < t \leq s$ we can conclude  $i \in S^{x'}$. In conclusion, since $i \in S^{x'}$ and due to~\eqref{eq:pi>=Pxyomega/Rd}, we also have $i \in S^{x''}$. Hence, $x''$ is indeed a seried strategy with $S^{x''} = \{1,\ldots,t\}$. This, together with~\eqref{eq:bijnaklaarstellingPxy>=Px''y''}, concludes the proof.
\end{proof} \bigskip

\newpage
\noindent \underline{\textbf{Proof of Theorem~\ref{algorithm}}}

\begin{proof}
First, note that it follows directly from the pseudocode that Algorithm~\ref{alg:algorithm}  is a linear time algorithm. Next, we show that Algorithm~\ref{alg:algorithm} finds an optimal strategy for the leader. For this, let $i^*$ be highest index for which $p_{i^*} > \overline{p}(\{1,2,\hdots,i^*-1\})$, i.e.,
$$i^* = \max\{j \in N \hspace{1mm} \vert \hspace{1mm} p_j > \overline{p}(\{1,2,\ldots,j-1\})\}.$$
Note that $i^*=i-1$, with $i$ the end value of variable $i$ in the algorithm. Moreover, note that $i^*$ always exists, since $p_1 > 0 = \overline{p}(\emptyset)$. The definition of $i^*$ implies that, for all $j \in \{1,2\ldots,i^*\}$, we have
$$p_j > \overline{p}(\{1,2,\ldots,j-1\}).$$
By combining this with Lemma~\ref{veiligtoevoegen}($i$), it follows that, for all $j \in \{1,2,\ldots,i^*\}$, we have
$$ \overline{p}(\{1,2,\ldots,j\}) \geq \overline{p}(\{1,2,\ldots,j-1\}) .$$ As a consequence,
\begin{equation} \label{Theorem2:hulp1}
\overline{p}(\{1,2,\ldots,i^*\})  \geq \overline{p}(\{1,2,\ldots,i^*-1\}) \geq \ldots \geq \overline{p}(\{1,2\}) \geq  \overline{p}(\{1\}) .
\end{equation}

We assume now, for the time being, that $i^*<n$ and thus index $i^*+1$ exists. By definition of $i^*$ and because the production rates are in non-increasing order,  we have
$$\overline{p}(\{1,2,\ldots,i^*\}) \geq p_{i^*+1}  \geq p_{i^*+2}\geq \ldots \geq p_{n-1}   \geq  p_n.$$
By combining this with Lemma~\ref{veiligtoevoegen}($ii$), it follows that, for all $j \in \{i^*+1,i^*+2,\ldots,n\}$, we  have
\begin{equation} \label{Theorem2:hulp2}  \overline{p}(\{1,2,\ldots, i^*\}) \geq \overline{p}(\{1,2,\ldots,j\} .\end{equation}

By combining \eqref{Theorem2:hulp1} and \eqref{Theorem2:hulp2}, we conclude that, for all $j \in \{1,2,\ldots,n\}$, we have \begin{equation}
\label{Theorem2:hulp3}
\overline{p}(\{1,2,\ldots,i^*) \geq \overline{p}(\{1,2,\ldots,j\}).
\end{equation}
Note that, if $i^*=n$, then~\eqref{Theorem2:hulp3} follows immediately from~\eqref{Theorem2:hulp1}.

Using Corollary~\ref{cor:allocationbalancedx+balancedNIEUW} together with~\eqref{Theorem2:hulp3}, we conclude that seried-balanced strategy $x^*$ with $S^{x^*}=\{1,2,\ldots,i^*\}$ has, among all seried-balanced strategy,  the highest worst case total production after destruction. From Lemma~\ref{bijnaklaarstelling} we know that there exists an optimal strategy for the leader that is seried-balanced. As a result, seried-balanced strategy $x^*$ with $S^{x^*} = \{1,2,\ldots,i^*\}$ is an optimal strategy for the leader. From  Lemma~\ref{allocationbalancedx+balanced} we know
\begin{equation*}
x^*_i = \begin{cases}
{\displaystyle \frac{a_iR_l}{p_i \sum_{j \in S^{x^*}} \frac{a_j}{p_j}}} & \text{ if } i \in S^{x^*}, \\
0 & \text{ if } i \in N \backslash S^{x^*},
 \end{cases}
\end{equation*}
which is in line with the output of Algorithm~\ref{alg:algorithm}.
\end{proof}

\newpage
\end{document}